\documentclass[reqno,11pt]{amsart}

\usepackage{tikz-cd}
\usepackage{tikz}

\usepackage{tcolorbox}

\usepackage{enumerate}
\usepackage{amsmath}
\usepackage{amssymb}
\usepackage{graphicx}
\usepackage[cmtip,all]{xy}
\usepackage[mathscr]{eucal}
\usepackage{bbold}            
\usepackage{tikz}
\usepackage{tikz-cd}

\newtheorem{theorem}{Theorem}[section]
\newtheorem{lemma}[theorem]{Lemma}
\newtheorem{corollary}[theorem]{Corollary}
\newtheorem{proposition}[theorem]{Proposition}
\newtheorem{properties}[theorem]{Properties}
\newtheorem{claim}{Claim}
\theoremstyle{definition}
\newtheorem{definition}[theorem]{Definition}
\newtheorem{example}[theorem]{Example}

\newtheorem{notation}[theorem]{Notation}

\theoremstyle{remark}
\newtheorem{remark}[theorem]{Remark}
\numberwithin{equation}{section}

\newcommand{\Lwedge}{\mathsf{\Lambda}}

\newcommand{\hatmu}{\widehat{\mu}}

\newcommand{\field}[1]{\ensuremath{\mathbb{#1}}}
\newcommand{\bba}{\field{A}}

\newcommand{\bbd}{\field{D}}
\newcommand{\bbe}{\field{E}}

\newcommand{\bbl}{\field{L}}

\newcommand{\bbp}{\field{P}}
\newcommand{\bbq}{\field{Q}}

\newcommand{\bbs}{\field{S}}

\newcommand{\bbu}{\field{U}}

\newcommand{\bbx}{\field{X}}

\newcommand{\bbz}{\field{Z}}

\newcommand{\uA}{\underline{A}}
\newcommand{\uB}{\underline{B}}

\newcommand{\place}[1]{\ensuremath{\mathbb{#1}}}

\newcommand{\bone}{\place{1}}
\newcommand{\calo}{\mathcal{O}}

\newcommand{\calz}{\mathcal{Z}}

\newcommand{\scrA}{\mathscr{A}}

\newcommand{\scrC}{\mathscr{C}}
\newcommand{\scrD}{\mathscr{D}}
\newcommand{\scrE}{\mathscr{E}}
\newcommand{\scrF}{\mathscr{F}}

\newcommand{\scrI}{\mathscr{I}}
\newcommand{\scrJ}{\mathscr{J}}

\newcommand{\scrL}{\mathscr{L}}

\newcommand{\scrO}{\mathscr{O}}

\newcommand{\scrQ}{\mathscr{Q}}

\newcommand{\scrU}{\mathscr{U}}

\newcommand{\scrZ}{\mathscr{Z}}

\newcommand{\fra}{\mathfrak{a}}
\newcommand{\frA}{\mathfrak{A}}
\newcommand{\frb}{\mathfrak{b}}

\newcommand{\frC}{\mathfrak{C}}

\newcommand{\frg}{\mathfrak{g}}

\newcommand{\frh}{\mathfrak{h}}

\newcommand{\frp}{\mathfrak{p}}

\newcommand{\frq}{\mathfrak{q}}

 \newcommand{\mba}{\mathbf{a}}
 \newcommand{\mbA}{\mathbf{A}}
 \newcommand{\mbb}{\mathbf{b}}
 \newcommand{\mbB}{\mathbf{B}}
 \newcommand{\mbc}{\mathbf{c}}
 
 \newcommand{\mbe}{\mathbf{e}}

 \newcommand{\mbh}{\mathbf{h}}

 \newcommand{\mbm}{\mathbf{m}}

 \newcommand{\mbp}{\mathbf{p}}

 \newcommand{\mbs}{\mathbf{s}}
 
 \newcommand{\mbt}{\mathbf{t}}
  \newcommand{\mbu}{\mathbf{u}}

 \newcommand{\mbv}{\mathbf{v}}

\newcommand{\mbW}{\mathbf{W}}
 \newcommand{\mbx}{\mathbf{x}}
 
\newcommand{\mby}{\mathbf{y}}

 \newcommand{\mbz}{\mathbf{z}}
 
 \newcommand{\mbzero}{\mathbf{0}}

\newcommand{\colim}[1]{\operatorname{colim}_{#1}}

 \newcommand{\ev}{\mathsf{ev}}
 \newcommand{\mat}[2]{\mathscr{M}_{#1 \times #2}}
 \newcommand{\la}{\langle}
 \newcommand{\ra}{\rangle}

\newcommand{\proj}[1]{\mathsf{Proj}\left( #1 \right)}
\newcommand{\detsch}[3]{\scrD^{#1}_{{#2} \to {#3}}} 

\newcommand{\fdetsch}[3]{{\hat\scrD}^{#1}_{{#2} \to {#3}}} 

\newcommand{\qbdlp}[2]{\scrQ^{\oplus #1}_{#2}} 

\newcommand{\bgl}[1]{B_{\bullet}GL_{#1}} 

\newcommand{\spec}[1]{\mathsf{Spec}\left(#1 \right)}

\newcommand{\osch}[1]{\mathsf{Sch}/{#1}}  

\newcommand{\sch}{\mathsf{Sch}}  
\newcommand{\sets}{\mathsf{Sets}} 
\newcommand{\Hom}[3]{\mathsf{Hom}_{#1}(#2, #3) } 
\newcommand{\inthom}[2]{\underline{\mathsf{Hom}}\left(#1, #2 \right)} 

\newcommand{\intdc}[3]{H^{#1}_{\scrD}(#3; \bbz(#2))}

\newcommand{\vtriang}{\Delta \hspace{-0.235cm}\mathord{\raisebox{-0.2\depth}{\scalebox{0.42}{\( \Delta\)} } }
\hspace{-0.26cm}\mathord{\raisebox{-0.8\depth}{\scalebox{0.64}{\( \Delta\)}}}} 

\newcommand{\simplex}[1]{\vtriang{\kern.22em}^{#1}} 

\newcommand{\hche}[3]{C^{#1}_{#3,#2}} 

\newcommand{\shcg}[3]{CH^{#1}(#3; #2)}  

\newcommand{\GL}[1]{GL_{#1}}
\newcommand{\glt}[2]{GL_{#1}^{\times #2}}

\newcommand{\sbcx}[3]{\mathcal{Z}^{#1}(#3;  #2)} 

\newcommand{\sbcxeq}[3]{\mathcal{Z}_{\text{eq}}^{#1}(#3; #2)}

\newcommand{\zhom}[3]{\bbz\mathsf{Hom}_{#1}(#2,#3)}

\makeatletter
\newcommand{\extp}{\@ifnextchar^\@extp{\@extp^{\,}}}
\def\@extp^#1{\mathop{\bigwedge\nolimits^{\!#1}}}
\makeatother

\newcommand{\Exterior}{\mathchoice{{\textstyle\bigwedge}}%
{{\bigwedge}}%
{{\textstyle\wedge}}%
{{\scriptstyle\wedge}}}

\newcommand{\hdiff}[1]{{}_{#1}{\hat \partial}}
\newcommand{\hcodiff}[1]{{}_{#1}{\hat \delta}}

\makeatletter
\newcommand*{\Relbarfill@}{\arrowfill@\Relbar\Relbar\Relbar}
\newcommand*{\xeq}[2][]{\ext@arrow 0055\Relbarfill@{#1}{#2}}
\makeatother

\title[Higher Chow Groups]{Explicit   Chern cycles in $BGL$ and $KV$-theory}

\author[Paulo Lima-Filho]{Paulo Lima-Filho}
\address{Department of Mathematics, Texas A\&M University, College Station, TX 77843}
\email{plfilho@math.tamu.edu}
\thanks{The author thanks the support and hospitality of the Isaac Newton Institute for Mathematical Sciences, during the program "K-theory, algebraic cycles and motivic homotopy theory" (KAH2, 2022), when this project came to fruition.}
\subjclass[2010]{14C35, 14F42, 19E15, 19E20}

\date{}
\dedicatory{}

\commby{}

\begin{document}

\begin{abstract}
Using   determinantal schemes, we construct explicit cycles in the higher Chow complex \( \sbcx{p}{*}{B_\bullet GL}\) that represent the universal Chern classes in the higher Chow groups \( \shcg{p}{0}{B_\bullet GL} \). As an application, we use these cycles, along with a canonical \emph{stable moving lemma} for Karoubi-Villamayor \(K\)-theory, to give a direct construction of the Chern class homomorphisms \( \mbc_{p,r} \colon KV_r(R) \to \shcg{p}{r}{\spec{R}}\) for regular algebras over a field \( \Bbbk\).

\end{abstract}

 \maketitle
\tableofcontents



\vfill\eject

\section*{Introduction}

In the foundational paper \cite{HG-RRHKT} H. Gillet developed  a general framework, on  suitable categories of schemes, to define a theory of Chern classes for higher algebraic  K-theory with values in any cohomology theory
satisfying natural conditions that hold in a broad class of cohomology theories. 
The universality of this approach established  \cite{HG-RRHKT} as a standard reference when cohomological methods are used to study algebraic \(  K \)-theory and various types of regulator maps; e.g.  \cite{MR704632}, \cite{MR899412}, \cite{MR944998}, \cite{MR1031903} and \cite{MR1621424}. 

This approach was used by S. Bloch \cite{Blo-HCG} to derive  the existence of universal Chern classes \(  C_p \in \shcg{p}{0}{\bgl{n}} \) with values in the Higher Chow groups of the classifying (simplicial) scheme \(  \bgl{n}.\)   The resulting Chern character, along with classical arguments, are then used to prove a Riemann-Roch isomorphism \(  \tau \colon G_n(X)_\bbq \xrightarrow{\cong} \oplus_i \shcg{i}{n}{X}_\bbq, \)  where \(  G_n(X) \) denotes the higher \(  K \)-theory of coherent sheaves on \(  X, \) promoting the higher Chow groups as the natural candidate for motivic cohomology in the category of   schemes over a field \( \Bbbk\).

The universal  classes, in turn, lead to  functorial Chern class homomorphisms \( \mbc_{p,r} \colon K_r(X) \to \shcg{p}{r}{X} \) in the category of regular schemes over \(\Bbbk\) (see Example \ref{ex:pairing1}) that are often coupled with realization morphisms \( \shcg{p}{r}{X} \to H^{2p-r}(X;A(p)) \) into suitable cohomology theories \( H^*(-;A(\bullet)) \) which are -- in principle -- amenable to computations. Studies of such realizations abound in the literature, amongst which Beilinson's regulator  \cite{MR862628} occupies a prominent status;  see  \cite{MR944987}, \cite{KL-AJHCG-II} and \cite{MR4498559}.
Nevertheless, a  simple and geometric presentation of the universal Chern classes in terms of concrete representing cycles seems to be missing from the literature, except for particular examples, such as  \cite{MR1237830}.

The primary goal of this paper is to provide such a construction, utilizing basic determinantal varieties, and  retrieving in a historical sense the classical intersection-theoretic origins of the subject; see \cite{Ful-IT}. Our ulterior motive was to use the Chern cycles to directly construct functorial Chern class homomorphisms 
\( \mbc_{p,r} \colon KV_r(R) \to \shcg{p}{r}{\spec{R}}\) from the Karoubi-Villamayor \( K\)-theory of a regular \( \Bbbk\)-algebra \( R\) to the higher Chow groups of \( \spec{R} \).  

The subtle difference in this approach is that we use the homotopy invariance of  \( KV_*(-)\)  together with the fact that an element \( \alpha \in KV_r(R) \) can be  {stably represented}, in a controlled way,  by a morphism with suitable {flatness properties}; see Corollary \ref{cor:KV-moving}. 
A key ingredient in the argument is \emph{Nori's equidimensionality lemma}, an unpublished result of Madhav Nori who kindly shared its proof and allowed us to include it as Appendix \ref{sec:LULU} in this paper.

In a nutshell, we prove a sort of \emph{moving lemma} for \( KV\)-theory that  holds for  integral \( \Bbbk\)-algebras of finite type. This allows for a direct pull-back of the Chern cycles constructed earlier.  This  is to be contrasted with utilizing the evaluation map
 \(
\ev \colon \mathsf{Spec}(R) \times B_\bullet GL_n(R) \longrightarrow B_\bullet GL_n
\)
 to pull-back the universal Chern classes in \( \shcg{*}{0}{B_\bullet GL}\), along with slant products and  Hurewicz maps. The needed contravariant functoriality here requires the regularity of \({R}\) along with a moving lemma for   higher Chow groups;  see Example  \ref{ex:pairing1}.

\medskip

For the reader's convenience and to establish notation, we provide a quick overview of higher Chow groups for simplicial schemes (with flat face maps) in Section \ref{sec:ssHCG}. To outline  our main constructions,  let  \(\mbx = (x_{ij}) \) be an \( n\times n\) matrix of indeterminates and let \( \mat{n}{n}:= \spec{\bbz[\mbx]}\) be the affine space  of \( n\times n \) matrices.  In Section \ref{sec:constr} we introduce and study a collection of morphisms
\(
L_{n,r} \colon \GL{n}^{\times r}\times \Delta^r \longrightarrow \mat{n}{n}
\) 
that send
\(  \mbA = (A_1|\cdots | A_r) \in \GL{n}(R)^{\times r} \)
and \( \mbt=(t_0, t_1, \ldots, t_r) \in \Delta^r(R)\) to
\begin{equation*}
L_{n,r}(\mbA, \mbt) :=  t_0 I + t_1 A_1 + t_2 (A_1A_2) +  \cdots + t_r(A_1\cdots A_r) \ \in \  \mat{n}{n}(R),
\end{equation*}
where \(  I = (\delta_{ij}) \) is the identity matrix, and
\(  
\Delta^r:= \spec{\bbz[z_0, \ldots, z_r]/\la z_0+ \cdots + z_r - 1 \ra}
\)
is the algebraic \(  r \)-simplex.
Denoting  the \( j\)-th column of  \( L_{n,r} \) by \( L_{n,r}^j\), 
we define  in Section \ref{sec:Ccycles} subschemes \(C^p_{n,r} \subset GL_n^{\times r} \times \Delta^r \), for each \( 1 \leq p \leq n \),   whose functor of points
 is 
\begin{equation*}
\label{eq:D-functor}
C^p_{n,r}(R) :=
 \left\{ 
 (\mbA,\mbt) \in GL_n(R)^{\times r} \times \Delta^r(R) \mid
 L_{n,r}^p(\mbA,\mbt)\wedge \cdots \wedge  L_{n,r}^n(\mbA,\mbt)  =  0 \in  \Exterior^{n-p+1}(R^{\oplus n}) 
\right\}.
\end{equation*}
It is worth mentioning that \( C^p_{n,r}\) is defined over \( \bbz\). The properties of this collection are summarized in Proposition \ref{prop:cycles2}, where we show that \( C^p_{n,r} \) is \emph{Cohen-Macaulay, normal, irreducible of codimension \( p\) in \( GL_n^{r}\times \Delta^r\) }. Furthermore, \emph{these subschemes are compatible under the inclusions \( GL_n\hookrightarrow GL_{n+1}\) and are equidimensional and dominant over \( \Delta^r\). }

The associated family of algebraic cycles \( \left\{ [C^p_{n,r}] \in \calz^p(GL_n^{\times r} \times \Delta^r) \mid r\geq 1\leq p \leq n\right\} \) satisfies  compatibility properties stemming from what we coin as \emph{a coherent family of cycles} in \( \mat{n}{n}\) in Section \ref{sec:constr}. 
These are conditions on families of algebraic cycles \( \{ \bbd^p_n \subset \mat{n}{n} \mid 1 \leq p \leq n\} \) that guarantee, amongst other things,  the existence of  flat  pull-backs \( L_{n,r}^*\, \bbd^p_n \)  and  that the collection \( ( L_{n,r}^*\, \bbd^p_n )_{r\geq 0}  \) defines a class in the higher Chow group \( \shcg{p}{0}{B_\bullet GL_n} \). 

We introduce  the \emph{Chern cycles} 
\(
\mathfrak{C}^p_n := ([C^{p}_{n,r}] ) \in \sbcxeq{p}{0}{B_\bullet GL_n} = \prod_r \sbcxeq{p}{r}{\glt{n}{r}}
\)
in Definition \ref{def:Chern_cyc} and prove the first main result.
\medskip

\noindent{\bf Theorem \ref{thm:main-T}.}{\it\ \ 
For each $1\leq p\leq n$, the element \( \mathfrak{C}^p_n \in \sbcxeq{p}{0}{\bgl{n}} \) is a cycle in the higher Chow complex of $\bgl{n}$ that represents the \(  p \)-th  Chern class \(  \mbc_p \in \shcg{p}{0}{\bgl{n}}  \) of the universal \(  n \)-plane bundle \(  \bbe_{n\, \bullet} \)  over \( \bgl{n}. \) Furthermore, the cycles \( \mathfrak{C}^p_n  \)  are compatible under  pull-back via the inclusions \(  \jmath_n \colon GL_n \hookrightarrow GL_{n+1} \).
}
\medskip

The proof of this theorem requires the study of certain families of determinantal schemes,  performed in Section \ref{sec:det_fam}, where we show that their associated cycles form coherent families in the sense of Section \ref{sec:constr}, and we determine how their associated cycles intersect.

In Section \ref{sec:CC-KV} we utilize the Chern cycles constructed above to directly present the aforementioned Chern class maps 
\( \mbc_{p,r} \colon KV_r(R) \to \shcg{p}{r}{\spec{R}} \) from \( KV\)-theory, whose construction goes roughly as follows. 
Let \( \bbl_{n}, \bbu_{n} \subset  SL_{n} \) be  the closed subgroups of \emph{lower triangular} and \emph{upper triangular unipotent} matrices, respectively, and denote $\bbx_n :=  \bbl_{n}\times \bbu_{n}\times \bbl_{n}\times \bbu_{n}$,  an affine space of dimension \( 2n(n-1) \). If  \( \mu  \colon \bbx_n \to SL_{n}$ is the multiplication map, it follows from \cite{SSV-Gauss} that $\mu$ is surjective, a fact that we combine with
the following equidimensionality result.
\medskip

\noindent{\bf Theorem \ref{thm:equid}.}{\, (M. Nori) \ \ \it
Let  \(  G \) be an algebraic group of dimension \(  n,   \) 
and let \(  f_i \colon X_i \to G \), \(  i=1,\dots, k \), be dominant morphisms from irreducible varieties \(  X_i \) to \(  G. \) Define \(  X := X_1 \times \cdots \times X_k \) and \(  F\colon X \to G \) by \(  F(x_1, \ldots, x_k) = f_1(x_1)\cdots f_k(x_k). \) If \(  k\geq n \), then \(  F \) is an equidimensional morphism. 
}
\medskip

It follows that whenever \( m\geq n^2-1\) the map
\begin{align*}
\mu_m \colon {\bbx_n \times \cdots \times \bbx_n} &\longrightarrow SL_{n} \\
\underline{\mbA}= (\mbA_1, \ldots, \mbA_m) & \longmapsto \mu(\mbA_1) \mu(\mbA_2) \cdots \mu(\mbA_m) \notag
\end{align*}
is \emph{faithfully flat}. Denoting  \(\scrF_n  := \underbrace{(\bbx_n \times \cdots \times \bbx_n)}_{n^2-1}  \ \cong \ \bba^{2n(n-1)(n^2-1)}\),  one can construct a ``contracting homotopy''  \( \frh_n \colon \scrF_n \times \bba^1 \to SL_n \)  whose properties are described next. 

\medskip

\noindent{\bf Proposition \ref{prop:H}.}{\ \ \it
 The following properties hold.
\begin{enumerate}[i.]
\item  \( \frh_n(\underline \mbA;0) = {I}\) for all \( \underline \mbA \in \scrF_n\), where \( {I}\) is the identity element in \( SL_n\).
\item  Given \( y \in \bba^1- \{ 0 \} \), the restriction \(  \frh_{n|y} \colon \scrF_{n|\Bbbk(y)} \to SL_{n|\Bbbk(y)}\) is faithfully flat.
\item   The natural inclusion \( \iota_n \colon \scrF_n \hookrightarrow \scrF_{n+1} \)  induced by \( \jmath_n \colon SL_n \hookrightarrow SL_{n+1} \), embeds \( \scrF_n \) as an affine subspace of \( \scrF_{n+1}\), yielding a commutative diagram.
\[
\begin{tikzcd}
\scrF_n \times \bba^1   \ar[rr, "\frh_n"] \ar[d, hook, "\iota_n \times \bone"'] &  & SL_n \ar[d, hook, "\jmath_n"] \\
\scrF_{n+1} \times \bba^1 \  \ar[rr, "\frh_{n+1}"] &  & SL_{n+1} 
\end{tikzcd}
\]
\end{enumerate}
}
\medskip

These properties along with other straightforward arguments yield the desired result, where \( \scrF_{n,r} := \scrF_{n}^{\times (r+1)} \) and \( X = \spec{R} \), with \( R \) a  \( \Bbbk\)-algebra.
\medskip

\noindent{\bf Corollary  \ref{cor:KV-moving}.}{\ \ \it
Let  \( \alpha \in KV_r(X), r \geq 2, \)  be represented by  \( \mba \colon X \times \Delta^r \to SL_n^{\times r} \). Then, there exists a canonical 
\( \lambda_\mba \colon X\times \scrF_{n,r} \times \Delta^r \to SL_n^{\times r} \) satisfying the following.
\begin{enumerate}[i.]
\item The restriction of \( \lambda_\mba\) to \( X\times \scrF_{n,r} \times (\Delta^r - \partial\Delta^r) \) is faithfully flat.
\item \( \lambda_\mba\) represents a class \( \lambda_\alpha \in KV_r(X\times \scrF_{n,r})\) corresponding to \( \alpha \) under the canonical isomorphism
\( \rho^*\colon KV_r(X) \xrightarrow{\ \cong \ } KV_r(X\times \scrF_{n,r}) \) induced by the projection \( \rho \colon X\times \scrF_{n,r} \to X \).
\end{enumerate}
}
\medskip

With the additional flatness property gained by this result, we  use \( \lambda_\mba\) to pull back the universal Chern cycle \( C^p_{n,r}\) to \( X\times \scrF_{n,r} \times \Delta^r \), thus defining \( \mbc_{p,r}\).
\medskip

\noindent{\bf Corollary  \ref{cor:chernKV}.}{\ \ \it 
For \( r \geq 2\) and \(  p\geq 1\), the  homomorphism
\(
\mbc_{p,r}  \)
is a natural transformation between the functors \(KV_r(-)\) and \(\shcg{p}{r}{\spec{-}} \) in the category of regular \(\Bbbk\)-algebras.
}
\medskip

When \( r=1\) we use separate considerations to deal with   \(K_1(R) \cong  R^\times \oplus SK_1(R) \). The factor \( SK_1(R) \)  is dealt with using similar flatness arguments, whereas for \( R^\times \) we resort to a rather amusing use of the assignment
\begin{equation*}
\alpha \in R^\times \mapsto A_\alpha := \begin{pmatrix} x t_0 & x^2 t_0t_1 - \alpha \\ 1 & xt_1 \end{pmatrix} \in GL_2(R[x][\Delta^1]),
\end{equation*}
to adhere with the use of Suslin's equidimensional complex to define higher Chow groups. 

We conclude the paper with a last application of the matrices \( A_\alpha \) and the universal Chern cycles, to construct algebraic cycles 
\([\Gamma_{\alpha, \beta} ] \) in \( X \times \bba^2 \times \Delta^2\) that yield the following.
\medskip

\noindent{\bf Corollary \ref{cor:symbol}.}{ \ \ \it
The assignment 
\[
  \alpha\otimes \beta  \in R^\times\otimes R^\times\ \overset{\gamma}{\longmapsto}\  [ \Gamma_{\alpha,\beta}] =: \gamma_{\alpha,\beta} \in \shcg{2}{2}{X\times \bba^2} \cong \shcg{2}{2}{X} 
\] 
is a Steinberg symbol defining a natural transformation 
\[
\gamma 
\colon K_2^M(R)  \to \shcg{2}{2}{\spec{R}} 
\]
of functors in the category of regular \( \Bbbk\)-algebras. 
}
\medskip

Although  constructions as in this corollary are  well-known, this direct approach  may have further interest.  Finally, we must point  out
the similarity between the definition of  \( L_{n,r} \) with  parametrizations of the geodesic simplices \( \Delta(\gamma_1, \ldots, \gamma_r) \subset G/K \)    used by J. Dupont  \cite{Dup-ccflat} to express the van-Est isomorphism  in terms of simplicial DeRham cohomology. One may na\"ively wonder if this underlies a DeRham realization of the geometric characteristic classes   presented here. 



\section{Preliminaries on simplicial schemes and Higher Chow Groups}
\label{sec:ssHCG}

We start with a brief digression about \emph{higher Chow groups for simplicial schemes} with flat face maps \cite[\S7]{Blo-HCG}, reviewing their functorial  properties, and describing exterior and intersection products in the simplicial context. 
When the simplicial scheme has the form $X\times Y_\bullet$, where $X$ is an arbitrary scheme and $Y_\bullet$ is a simplicial set, we discuss the \textit{slant product} pairing from the higher Chow groups of $X\times Y_\bullet$ and the singular homology of $Y_\bullet$ ($\bbz$ coefficients) to the higher Chow groups of $X.$  

\subsection{Conventions and notation}

Throughout this paper  \( \Bbbk\) is a field and, \emph{unless otherwise stated}, all schemes  have locally finite type and are defined over \(  \mathbb{k}. \) We let $\osch{\mathbb{k}}$ denote the category of such schemes and morphisms over \(\Bbbk\).  Subscripts in the notation of categorical products are omitted, except when needed for clarity. Whenever convenient, we use Yoneda's identification between an algebraic scheme \( X\) and its functor of points \( R\mapsto X(R):= \Hom{\osch{\Bbbk}}{\spec{R}}{X} \) on the category of \( \Bbbk\)-algebras. 

The symbol \( \Exterior^k S \) will    either   denote the \( k\)-th exterior power of an \(R\)-module \( S\), or the subset of the power set of \( S\) consisting of the finite subsets  with \( k \) elements. Most notably, the latter meaning is used when \( S = [n] := \{ 1, \ldots, n\}\).

We refer the reader to \cite[\S 1]{Ful-IT} for basic terminology and results on algebraic cycles, from which we borrow the following notation.  If  \( Z \) is a closed \( d\)-dimensional subscheme of a scheme \( X\), then \( [Z]\) denotes the corresponding algebraic \(d\)-cycle on \(X \). If the codimension of an algebraic cycle in \(X\) is given, we implicitly assume that \(X\) is an equidimensional scheme. 

Let \( f \colon X \to Y \) be a morphism of schemes and let \( \mba \) and \( \mbb\)  be algebraic cycles  in \(X \) and \(Y\), respectively. If \( f\) is a proper morphism, denote by \( f_*\mba\) the proper push-forward of \( \mba\) under \( f \). Similarly, if \( f \) is flat of relative dimension \( e \),   denote by \( f^*\mbb \) the flat pull-back of \( \mbb\) under \( f\). We use the same notation to denote the corresponding operations on (higher) Chow groups, whenever applicable.  If \( V \subset Y\) is a closed subscheme,   \( f^{-1}Y \) is the scheme-theoretic inverse image of \( Y\). We also use \( f^* \)  to denote pull-back of equidimensional cycles, in the sense of \cite{Sus&Voe-RelCh}, and the context will make its meaning clear. 

If the algebraic cycles \( \mbu, \mbv\) meet properly in \(X\), then  \( \mbu \centerdot \mbv\)  denotes their \emph{intersection-theoretic} intersection,  while the symbol \( \cap \)  is reserved for set-theoretic or scheme-theoretic intersections. On regular schemes,  \( \alpha \cup \beta\) is the product of elements \( \alpha, \beta\) in the (higher) Chow ring \( \shcg{*}{\bullet}{X}\).

\subsection{Simplicial schemes with flat face maps}

We first recall the definition of higher Chow groups for simplicial schemes with flat face maps, taken essentially from \cite[\S7]{Blo-HCG}. The only departure from the original formulation lies in the use of Suslin's equidimensional complexes, along with  products instead of direct sums. 

Denote the co-faces and  co-degeneracies of $\Delta^\bullet$ by
\begin{equation}
\label{eq:co-simp}
\partial_j \colon \Delta^{n-1} \to \Delta^n \quad \text{and} \quad s_j \colon \Delta^{n+1} \to \Delta^{n}, \ \ j = 0, \ldots, n,
\end{equation}
respectively.  Consider a simplicial scheme \( X_\bullet  \) whose face maps $\delta_j \colon X_n \to X_{n-1}$, $j=0, \ldots n, $ are flat morphisms. For a fixed \(  p\geq 0 \), the flat pull-back of cycles yields a double complex with groups 
$\sbcxeq{p}{s}{X_r}
$
and commuting differentials
\[
\xymatrix{
\sbcxeq{p}{s}{X_r} \ar[r]^-{\hcodiff{r}} \ar[d]_{\hdiff{s}} & \sbcxeq{p}{s}{X_{r+1}} \ar[d]^{\hdiff{s}} \\
\sbcxeq{p}{s-1}{X_r} \ar[r]_-{\hcodiff{r}}
 & 
\sbcxeq{p}{s-1}{X_{r+1}},
}
\]
where  $\hdiff{s} := \sum_{i=0}^s (-1)^i ( \bone \times \partial_i)^*$  are the differentials in the higher Chow complexes and $\hcodiff{r} := \sum_{j=0}^{r+1}(-1)^j (\delta_j \times \bone)^* $ comes from the face maps of the simplicial scheme. Subscripts will be omitted when the context is clear. 

\begin{definition}
\label{def:HCGssch}
Let $X_\bullet$ be a simplicial scheme with \emph{flat face maps}. The higher Chow complex $\sbcxeq{p}{*}{X_\bullet}$
is the total (product) complex   associated to the double complex above. In other words, 
\[ 
\sbcxeq{p}{k}{X_\bullet}
:= \prod_{s-r=k} \sbcxeq{p}{s}{X_r} = \prod_{r\geq 0} \sbcxeq{p}{r+k}{X_r},
\]
and its differential 
\(
d\ \colon \  \sbcxeq{p}{k}{X_\bullet}  \longrightarrow \sbcxeq{p}{k-1}{X_\bullet}
\)
takes $\beta = (\beta_r)$\  to\  $d\beta =(\{ d\beta\}_{i})$, where $\{d\beta\}_{i}:= \hat \partial \beta_{i} \ + \ (-1)^{k-1}\, \hat\delta \beta_{i-1}$. The $k$-th higher Chow group of $X_\bullet$\  of weight $p$ is the $k$-th homology group 
\(
CH^p(X_\bullet;k) := H_k(\sbcxeq{p}{*}{X_\bullet})
\).

\end{definition}
\begin{example}
If $Y_\bullet$ is a simplicial set and $X$ is a scheme of finite type over $\Bbbk$, then $X \times Y_\bullet$ becomes a simplicial scheme in $\osch{\Bbbk}$,  whose face maps are flat.
\end{example}
\begin{example}
\label{exmp:bar}
Let $G$ be an algebraic group over a base scheme $S$. A \emph{flat $G$-triple $( X, G, Y) $ over $S$} consists of schemes  $X, Y$ in $\osch{S}$ endowed with a right action of $G$  on $X$ and a left action on \(  Y, \) having the property that the action morphisms
$\rho \colon X\times G \to X$ and $\lambda \colon G \times Y \to Y$ are flat.
The \emph{triple bar construction} yields a simplicial scheme $B_\bullet(X,G,Y)$ with  $B_r(X,G,Y):= X\times_S G \times_S \cdots \times _S G \times_S Y$, whose face maps are flat. See  \cite{May-classifying} or \cite[\S 2.3]{LF-TPACF}. Particular examples include:

\begin{align*}
B_\bullet G_{/S} & := B_\bullet(S,G,S) = \text{ classifying space of } G_{\!/S} \\
\bbe_{n\bullet{/S}}& := B_\bullet(\bba^n_S, GL_{n/S},S) \to B_\bullet GL_{n/S}  \\
Gr^k_{n\, \bullet/S} & := B_\bullet(Gr_{k,n/S}, GL_{n/S},S) \to B_\bullet GL_{n/S}.
\end{align*}
In other words, $\bbe_{n,\bullet/S}$ and $Gr^k_{n\, \bullet/S}$ \ are, respectively, the \emph{universal \(n\)-plane bundle} and the \emph{Grassman bundle} of codimension $k$ quotients of  $\bbe_{n,\bullet}$ over $B_\bullet GL_{n/S}$.
\end{example}

\subsection{Products and functoriality}

This section translates basic combinatorial constructions, such as the Eilenberg-Zilber and Alexander-Whitney maps and exterior products, into the context of higher Chow groups. Working at the level of   equidimensional complexes makes this translation a mere formality that we pursue only to introduce the  notation and lemmata needed subsequently. The sole additional ingredient is   the the notion of  \emph{special cycles} (for the lack of better terminology), that we show to behave particularly well with respect to products.

Given a category \(  \scrC  \),  let   \(\,  \bbz\scrC\, \) denote the free \(  \mathsf{Ab}\)-category generated by $\scrC,$ having
\( \mathsf{Ob}(\bbz\scrC) = \mathsf{Ob}(\scrC)\)  as objects and   
\( \Hom{\bbz\scrC}{X}{Y} := \zhom{\scrC}{X}{Y} \)  as morphisms, 
where \( \bbz S \) is the free abelian group on  a set \(S. \) Composition gives pairings
 \begin{equation}
     \label{eq:pair}
\circ \colon   \Hom{\bbz\scrC}{Y}{Z}\otimes \Hom{\bbz\scrC}{X}{Y} \longrightarrow \Hom{\bbz\scrC}{X}{Z}.
\end{equation}

\subsubsection{Combinatorial jollies}

For the reader's convenience we recall classical combinatorial constructions going back to \cite{EZ-products}, working  directly in the category $\bbz\sch.$

\begin{definition}
\label{def:AWEZ} 
Consider  non-negative integers $a, b, \ell$ satisfying $a+b=\ell.$ 
\begin{enumerate}[a.]
\item \label{it:EZ} The \emph{Eilenberg-Zilber} map 
$\psi_{a,b} \in \Hom{\bbz \sch }{\Delta^\ell}{\Delta^a \times \Delta^b}$
is defined as 
\begin{equation}
    \label{eq:EZ}
    \psi_{a,b} := \sum_{I \in \Exterior^a[\ell-1]} \varepsilon_I\, (\mbs_{I^\complement} \times \mbs_{I})
\end{equation}
where   $\mbs_J := s_{j_1} \cdots s_{j_k} \colon \Delta^\ell \to \Delta^{\ell-k}$ 
for $J = (j_1< \cdots < j_k) \subset [\ell-1] $,  
and $\varepsilon_I$ is the sign of the shuffle permutation for the partition $[\ell-1] = I \amalg I^\complement $. Note that $\mbs_{I^\complement} \times \mbs_{I} \colon \Delta^\ell \to \Delta^a \times \Delta^b$ is an isomorphism in $\sch.$
\item \label{it:AW} The Alexander-Whitney map $E_{a,b} \in \Hom{\bbz \sch}{\Delta^a \times \Delta^b}{\Delta^\ell\times \Delta^\ell}$ 
is defined as 
\begin{equation}
    \label{eq:AW}
E_{a,b} := (\partial_{\ell}\cdots \partial_{a+1}) \times \partial_0^a \ = \ (\partial_{\ell}\cdots \partial_{a+1}) \times (\partial_0  \cdots \partial_{a-1} )
\end{equation}
with $\partial_{\ell}\cdots \partial_{a+1} \in \Hom{}{\Delta^a}{\Delta^\ell}$ and $\partial_0^a = \partial_0 \cdots \partial_0 = \partial_0  \cdots \partial_{a-1} \in \Hom{}{\Delta^b}{\Delta^\ell}.$
\item \label{it:nabla} Define $\nabla_\ell \colon \Delta^{\ell-1} \times \Delta^{\ell-1} \to \Delta^{\ell} \times \Delta^{\ell}$ by
$\nabla_\ell := \sum_{j=0}^{\ell+1} (-1)^j \partial_j \times \partial_j$. This is the differential of the Moore complex associated to the cosimplicial abelian group $\bbz(\Delta^\bullet \times \Delta^\bullet).$
\item \label{it:diag_aprox} The \emph{diagonal approximation} maps $\Phi_\ell \colon \Delta^\ell \to \Delta^\ell \times \Delta^\ell$ in $\bbz\sch$
are defined by
\begin{equation}
    \label{eq:d_approx}
\Phi_\ell := \sum_{a+b=\ell} E_{a,b} \circ \psi_{a,b}.
\end{equation}
\item \label{it:htpy} Finally, define $H_\ell \in \Hom{\bbz\sch }{\Delta^1\times \Delta^\ell}{\Delta^{\ell}\times \Delta^\ell}$ as the unique (linear combination of) affine linear map(s) that sends an $R$-valued point  $(\mbe_0,\mbx) \in \Delta^1(R) \times \Delta^\ell(R)$ to $\text{diag}_\ell(\mbx) = (\mbx, \mbx)$ and sends $(\mbe_1, \mbx)$ to $\Phi_\ell(\mbx).$

\end{enumerate}
\end{definition}
\begin{properties}
\label{proper:AWEZ}
Let $a, b, \ell$ be non-negative integers satisfying $a+b=\ell.$ The following identities hold. 
\begin{align}
    \psi_{a,b}\circ \hdiff{\ell} & = 
    (\hdiff{a}\times \bone)\circ \psi_{a-1,b} + (-1)^a (\bone\times \hdiff{b})\circ \psi_{a,b-1}, \label{eq:EZ1} \\
    \nabla_{\ell+1} \circ E_{a,b} & = E_{a+1,b} \circ (\hdiff{a+1} \times \bone) + (-1)^a E_{a,b+1}\circ (\bone \times \hdiff{b+1})
\end{align}
These relations are illustrated in the diagram below. 
\smallskip
\[
\xymatrix{
& 
\Delta^{a-1}\times \Delta^{b}   
    \ar[dr]^-{\hdiff{a}\times \bone} 
      \ar@{}[d] | {+}  
& & 
\Delta^{a+1}\times \Delta^{b} 
    \ar[dr]^{E_{a+1,b}}
    \ar@{}[d] | {+} 
& \\
\Delta^{\ell-1} 
    \ar[ur]^-{\psi_{a-1,b}}
    \ar[dr]_-{\psi_{a,b-1}}
    \ar[r]^{\hdiff{\ell}}
&
\Delta^{\ell} 
    \ar[r]^-{\psi_{a,b}} 
&
\Delta^{a}\times \Delta^b 
    \ar[r]^{E_{a,b}}
    \ar[ur]^{{}_{a+1}{\hat \partial}\times \bone }
    \ar[dr]_{\bone \times {}_{b+1}{\hat \partial} \ }
&
\Delta^{\ell}\times \Delta^\ell 
    \ar[r]^-{\nabla_{\ell+1}} 
    \ar@{}[d] | {+(-1)^a} 
&
\Delta^{\ell+1}\times \Delta^{\ell+1} 
\\
& 
\Delta^{a}\times \Delta^{b-1} 
    \ar[ur]_-{\bone \times {}_b{\hat \partial}}
       \ar@{}[u] | {+(-1)^a}   
& & 
\Delta^{a}\times \Delta^{b+1} 
  \ar[ur]_-{E_{a,b+1}}
& 
}
\]
\end{properties}
\medskip

\begin{lemma}
\label{lem:diag_approx}
The maps $\Phi_\ell$ \eqref{eq:d_approx} give a morphism of Moore complexes $\Phi_* \colon \bbz \Delta^* \to \bbz(\Delta^* \times \Delta^*).$ Hence, for all $\ell \geq 1,$ the diagram below commutes
\begin{equation}
    \label{eq:diag_cx}
    \xymatrix{
    \Delta^\ell \ar[rr]^{\Phi_\ell} & & \Delta^\ell\times \Delta^\ell \\
    \Delta^{\ell-1} \ar[rr]_{\Phi_{\ell-1}} \ar[u]^{\hdiff{\ell}} & & \Delta^{\ell-1} \times \Delta^{\ell-1} \ar[u]_{\nabla_\ell}.
    }
\end{equation}
Now, if $P_{k+1} \in \Hom{\bbz\sch }{\Delta^{k+1}}{\Delta^k \times \Delta^k}$ is the composition $P_{k+1}:= H_{k+1} \circ \Psi_{1,k}, $
then the following formula holds for all $k\geq 0.$
\begin{equation}
    \label{eq:htpy}
    \Phi_\ell  -  \operatorname{diag}_\ell \ = \ P_{\ell+1}\circ \hdiff{\ell+1}\, +\,  \nabla_\ell\circ P_{\ell},
\end{equation}
i.e., $P_*$ gives a homotopy between $\Phi_*$ and $\operatorname{diag}_*$ at the level of Moore complexes. 
\end{lemma}

\begin{remark}
\label{rem:prism}
Later, we will need the following repackaging  of the classic prism construction .
Suppose that  \( X_\bullet \)  is a simplicial scheme with flat face maps and fix \( k \geq 0 \). For each \( r \geq 0 \) let  \( \Gamma_r \) be an algebraic cycle in  \(  X_r \times \Delta^{r+k}\times \bba^1 \)  which  has codimension \( p \) and is equidimensional and dominant over \( \Delta^{r+k} \times \bba^1 \).

Denote  \( A_r := \Gamma_{r|_{T=1}} \) and \( B_r := \Gamma_{r|_{T=0}} \), so that \( \mbA := (A_r ) \) and \( \mbB:= (B_r) \) are elements in \( \sbcxeq{p}{k}{X_\bullet}\). Identifying \(\bba^1 \) with \( \Delta^1 \) via \( T \mapsto (1-T,T)\), we use the Eilenberg-Zilber maps to define
\( W_r := (-1)^{r+k}(\bone \times \psi_{r+k,1})^*\Gamma_r \), and   \( \mbW := (W_r) \in \sbcxeq{p}{k+1}{X_\bullet} \). 
%

Now, for each \( x \in \bba^1\), let \( \Gamma_{r |\Bbbk(x)} := \Gamma \centerdot \{  X \times \Delta^r \times \spec{\Bbbk(x)} \} \) be the ``slice'' of \(\Gamma_r\) above \( x\). If \( \Gamma_{|\Bbbk(x)} := ( \Gamma_{r|\Bbbk(x)} ) \) is a cycle in 
\( \sbcx{p}{k}{X_{\bullet | \Bbbk(x)}} \) for all \( x \in \bba^1\) , then it follows
from    \eqref{eq:EZ1} and Definition \ref{def:HCGssch} that
\[
d \mbW \ = \ \mbA  -  \mbB.
\]
\end{remark}

\subsubsection{Products in Higher Chow Groups of simplicial schemes}
The equidimensional higher Chow complexes $\sbcxeq{p}{*}{-}$ have the added benefit that all combinatorial maps discussed above induce   morphisms and homotopies in the level of complexes. This feature is explored in the discussion of products and \textit{special cycles} in the higher Chow complexes of simplicial schemes.

\begin{definition}
\label{def:prod_simp}
Let $X_\bullet, Y_\bullet $ be simplicial schemes with flat face maps 
$$\delta_i \colon X_\ell \to X_{\ell-1} \quad\text{ and } \quad  \delta'_j \colon Y_\ell \to Y_{\ell-1}, \ \ j = 0, \ldots, \ell-1, $$
respectively. For $a+b = \ell$, let $S_{a,b} \colon X_\ell \times Y_\ell \to X_a \times Y_b$ be the flat morphism 
 $S_{a,b} := (\delta_a \cdots \delta_{\ell}) \times \delta'^{a}_0$. 
Given $\alpha = (\alpha_r) \in \sbcxeq{p}{m}{X_\bullet }$ and $\beta= (\beta_s) \in \sbcxeq{q}{n}{Y_\bullet}$, define their \emph{exterior product}
$\alpha\boxtimes\beta    \in \sbcxeq{p+q}{m+n}{X_\bullet \times Y_\bullet}$ by
\begin{equation}
    \label{eq:boxprod}
    \{ \alpha\boxtimes\beta\}_\ell := \ \sum_{a+b=\ell} (-1)^{an}\, \left( S_{a,b} \times \psi_{a+m,b+n} \right)^*(\alpha_a \times \beta_b).
\end{equation}
Here $\alpha_a \times \beta_b $ is the product of cycles, under the identification 
\[
X_a\times Y_b \times \Delta^{a+m}\times \Delta^{n+b} \equiv (X_a\times \Delta^{a+m}) \times ( Y_b \times \Delta^{n+b}),
\] and   $S_{a,b} \times \psi_{a+m,b+n}$ is seen as a morphism in 
$\bbz \osch{\Bbbk}$ from $X_\ell\times Y_\ell \times \Delta^{\ell+m+n}$\  to\  
$X_a\times Y_b \times \Delta^{a+m}\times \Delta^{n+b}$. The pull-back maps \( ( S_{a,b} \times \psi_{a+m,b+n} )^* \) are well-defined, since the  face maps are flat, and the cycles used are suitably equidimensional.
\end{definition}
\begin{proposition}
\label{prop:ext_pair}
Given $X_\bullet$ and $Y_\bullet$ as in the definition, the pairing
\[
\boxtimes \colon \sbcxeq{p}{*}{X_\bullet} \otimes\sbcxeq{q}{*}{Y_\bullet} \longrightarrow \sbcxeq{p+q}{*}{X_\bullet\times Y_\bullet}, \quad \alpha\otimes \beta \mapsto \alpha \boxtimes\beta,
\]
is a map of complexes, i.e., for $\alpha \in \sbcxeq{p}{m}{X_\bullet}, \beta  \in \sbcxeq{n}{n}{X_\bullet}$ one has 
$$ d(\alpha\boxtimes \beta) = (d\alpha)\boxtimes \beta + (-1)^m \alpha \boxtimes (d\beta).$$
Furthermore, this exterior product is strictly associative. 
\end{proposition}

Using Suslin's generic equidimensionality results \cite{Sus-HighCh} and a simple argument with spectral sequences  one can also prove the following result.
\begin{proposition}
\label{prop:moving}
Let $j\colon V_\bullet \hookrightarrow X_\bullet$ be a regularly embedded simplicial subvariety of codimension $c$, where both $V_\bullet$ and $X_\bullet$ have flat face maps. Denote by $\sbcxeq{p}{*}{X_\bullet}_{V_\bullet}$ the subcomplex of  $\sbcxeq{p}{*}{X_\bullet}$
generated by  those cycles that intersect $V_\bullet$ properly. Then the pull-back map induced by the classical intersection product
\begin{equation}
    \label{eq:pull-back}
    j^* \colon \sbcxeq{p}{*}{X_\bullet}_{V_\bullet} \longrightarrow \sbcxeq{p}{*}{V_\bullet}, \quad \alpha \mapsto \alpha \centerdot  [Y_\bullet] := \left( \alpha_r  \centerdot  [Y_r \times \Delta^{*+r}]\right),
\end{equation}
is a quasi-isomorphism.
\end{proposition}

\begin{remark}
\label{rem:int_funct}
It follows from Proposition \ref{prop:moving}   that the higher Chow groups form  a contravariant functor in the category of regular simplicial schemes with flat face maps. Furthermore, if
$\operatorname{diag}  \colon X_\bullet \to X_\bullet \times X_\bullet$ is the diagonal map,  the assignment $[\alpha]\otimes [\beta] \mapsto [\alpha] \cup [\beta] := \operatorname{diag}^*\left[ \alpha\boxtimes \beta \right]$
defines an associative intersection product
\begin{equation}
    \label{eq:simpl_cup}
    \cup \colon \shcg{p}{m}{X_\bullet}\otimes \shcg{q}{n}{X_\bullet} \longrightarrow \shcg{p+q}{m+n}{X_\bullet}
\end{equation}
satisfying $\alpha\cup \beta = (-1)^{mn}\beta\cup \alpha. $ It follows that  $\shcg{*}{*}{X_\bullet}$ has the structure of a (bi)graded-commutative ring.
\end{remark}

\subsubsection{Products of special cycles}
Next, we introduce the notion of \emph{special cycles},  which  are elements in \( \sbcxeq{k}{0}{V_\bullet} \)  for a given simplicial scheme \( V_\bullet \), satisfying simple conditions with respect to the boundaries of the original double complex. These conditions imply that  exterior products of  classes represented by such cycles in the higher Chow groups have particularly nice representatives.

Given \(  \alpha = (\alpha_r) \in \sbcxeq{p}{0}{X_\bullet} \) and \(  \beta = (\beta_s) \in \sbcxeq{q}{0}{Y_\bullet} \), let
\begin{equation}
\label{eq:star_prod}
(\alpha \star \beta)_r := (\bone_{X_\bullet} \times \bone_{Y_\bullet} \times \operatorname{diag}^\Delta_r)^*(\alpha_r \times \beta_r) \in \sbcxeq{p+q}{r}{X_r \times Y_r}
\end{equation}
be the components of an element $\alpha \star \beta \in \sbcxeq{p+q}{0}{X_\bullet \times_\Bbbk Y_\bullet}$.
\begin{definition}
\label{def:special_cycle}
Let $V_\bullet$ be a simplicial scheme. An element \(  \gamma \in \sbcxeq{k}{0}{V_\bullet} \) 
is called a \emph{special cycle} if
\begin{equation}
\label{eq:special}
(\bone_{V_r} \times \partial_j)^* \gamma_r  = ( \delta_j  \times \bone_{\Delta^{r-1}} )^*(\gamma_{r-1}),
\end{equation}
for all \( 0 \leq j \leq r .\) 
Note that this condition implies that \(  d \gamma = 0 \).
\end{definition}
\begin{proposition}
\label{prop:special}
Given special cycles \(  \alpha \in \sbcxeq{p}{0}{X_\bullet} \) and \(  \beta \in \sbcxeq{q}{0}{Y_\bullet} \) 
then \(  \alpha \star \beta \) is also a cycle, and it represents the exterior product of their homology classes. In other words,  
\( [\alpha]\boxtimes [\beta] := [ \alpha \boxtimes \beta]  =  [\alpha \star \beta] \in CH^{p+q}(X_\bullet\times_\Bbbk Y_\bullet; 0)\).
\end{proposition}
\begin{proof}
By definition, 
\begin{align*}
    (\bone_{X_r} & \times \bone_{Y_r} \times \partial_j)^*(\alpha\star\beta)_r  = 
     (\bone \times\bone \times \partial_j)^*(\bone\times\bone \times \operatorname{diag}_r^\Delta)^*(\alpha_r
    \times \beta_r) \\
    & = (\bone   \times \partial_j)^*\alpha_r \times (\bone  \times \partial_j)^*\beta_r  
     = (\delta_j   \times \bone)^*\alpha_{r-1} \times (\delta'_j  \times \bone)^*\beta_{r-1} \\
    & = (\delta_j \times \delta'_j \times  \bone)^* (\bone\times \bone \times \operatorname{diag}^\Delta_r)^*(\alpha_{r-1} \times \beta_{r-1}) =  (\delta_j \times \delta'_j \times  \bone)^*\{ \alpha \star \beta\}_{r-1},
\end{align*}
and it follows that
\(
    \{d(\alpha\star \beta)\}_r =  \hdiff{} ( \alpha\star\beta)_r - \hcodiff{}(\alpha\star \beta)_{r-1} = 0.
\)

Using induction and Definition \ref{def:special_cycle} one  shows that whenever $a+b = \ell,$
the identity $(S_{a,b}\times \bone)^*(\alpha_a \times \beta_b) = (\bone \times E_{a,b})^*(\alpha_\ell\times \beta_\ell)$ holds. Therefore,
\begin{align}
    \{\alpha & \boxtimes\beta\}_\ell  =  \sum_{a+b=\ell} (S_{a,b}\times \psi_{a,b})^*(\alpha_a\times \beta_b) = \sum_{a+b=\ell} ( \bone \times \psi_{a,b})^*\circ  (S_{a,b}\times \bone)^*(\alpha_a\times \beta_b) \notag \\
   &= \sum_{a+b=\ell} ( \bone \times \psi_{a,b})^*\circ  (\bone \times E_{a,b})^*(\alpha_\ell\times \beta_\ell) = 
   \left(  \sum_{a+b=\ell}  \bone \times E_{a,b}\circ \psi_{a,b}\right)^*(\alpha_\ell\times \beta_\ell) \notag \\
\label{eq:boxt1}
& =
 (\bone\times \Phi_\ell)^*(\alpha_\ell\times \beta_\ell),
\end{align}
where $\Phi_\ell$ is the diagonal approximation map; see Definition \ref{def:AWEZ}\eqref{it:diag_aprox}.

Let $\theta = (\theta_\ell) \in \sbcxeq{p+q}{1}{X_\bullet\times Y_\bullet}$ be the element whose components are defined as $\theta_\ell := (\bone \times \bone \times P_{\ell+1})^*(\alpha_\ell\times \beta_\ell) \in \sbcxeq{p+q}{1}{X_\ell\times Y_\ell}$, where $P_{\ell+1} \colon \Delta^{\ell+1} \to \Delta^\ell \times \Delta^\ell$ is defined in Lemma \ref{lem:diag_approx}. 
It follows from \eqref{eq:htpy}, \eqref{eq:boxt1} and definitions that 
\begin{align*}
\{ d\theta\}_\ell &= (\bone \times \Phi_\ell)^*(\alpha_\ell\times \beta_\ell) - (\bone \times \operatorname{diag}^\Delta_\ell)^*(\alpha_\ell\times \beta_\ell)  = \{ \alpha\boxtimes \beta\}_\ell - \{ \alpha\star \beta\}_\ell,
\end{align*}
in other words, $d\theta = \alpha\boxtimes \beta \ - \ \alpha\star \beta$. 
\end{proof}

\subsection{Higher Chow groups and slant products}
In this section we consider simplicial schemes of the type \( X \times Y_\bullet\), where \( X\) is a scheme of finite type and \( Y_\bullet\) is a simplicial set. The main goal is to discuss the slant product from the higher Chow groups of  \( X \times Y_\bullet\) and the singular homology of \( Y_\bullet \) -- with \( \bbz \) coefficients -- to the higher Chow groups of \( X\).

Given a simplicial set $Y_\bullet$ and let $(\bbz Y_*, \mathsf{d} ) $ be the Moore complex of the simplicial group $\bbz Y_\bullet,$   which computes the integral singular homology of $Y_\bullet$.
If $X$ is a scheme of finite type over $\mathbb{k}$, for each $r\geq 0$ the group
$\sbcxeq{p}{k}{X\times Y_r}$ is a direct sum
$\sbcxeq{p}{k}{X\times Y_r}= \bigoplus_{y \in Y_r}\, \sbcxeq{p}{k}{X\times y\,}$, where $X\times y$ denotes a copy of $X$ indexed by $y \in Y_r$. 
Hence, one can write $Z \in \sbcxeq{p}{k}{X\times Y_r}$ as a sum 
$Z = \sum_{y\in Y_r} Z_y,$ with finitely many non-zero components $Z_y \in \sbcxeq{p}{k}{X\times y\,}$.

Let $\pi_y \colon X\times y \xrightarrow{} X$ be the  natural identification and 
define a pairing
\begin{equation}
\label{eq:slant}
\text{-} /\text{-} \ \colon\  \sbcxeq{p}{k}{X\times Y_r} \otimes \bbz Y_r \longrightarrow \sbcxeq{p}{k}{X}.
\end{equation}
by sending $Z = \sum_{y} Z_y\  \in\  \sbcxeq{p}{k}{X\times Y_r}$\ and \ $c= \sum_{g} n_g\cdot g \ \in \ \bbz Y_r$ to
\[
Z/c := \sum_{y,g \in Y_r}\ \delta_{y,g} \ n_g\ \pi_{y *}(Z_y) \ = \ \sum_{y \in Y_r} n_y\, \pi_{y *}( Z_y) \in \ \sbcxeq{p}{k}{X}
\]
where $\delta_{y,g} $ is the Kronecker delta.

Now, let 
\(
\rho_r \colon \sbcxeq{p}{k}{X\times Y_\bullet} = \prod_{i} \sbcxeq{p}{k+i}{X\times Y_i}  \longrightarrow \sbcxeq{p}{k+r}{X\times Y_r}
\)
be the projection on the \( r\)-th factor. Using \eqref{eq:slant} one obtains a pairing
\begin{equation}
\label{eq:slant_s}
    \text{-}/\text{-} \ \colon\  \sbcxeq{p}{k}{X\times Y_\bullet} \ \otimes \ \bbz Y_r \longrightarrow 
    \sbcxeq{p}{k+r}{X}
\end{equation}
defined by $\alpha \otimes c \mapsto  (\rho_r \alpha)/c.$

\begin{proposition}
\label{prop:slant}
The pairing above induces a slant product of complexes
\begin{equation}
    \label{eq:pairing}
    \text{-}/\text{-} \ \colon \ \sbcxeq{p}{*}{X\times Y_\bullet} \ \otimes \ \bbz Y_* \longrightarrow 
    \sbcxeq{p}{*}{X}
\end{equation}
which, in turn,  gives the \emph{slant product} in the level of homology  groups:
\begin{equation}
    \label{eq:slant_hom}
    \text{-}/\text{-} \ \colon \ CH^p(X\times Y_\bullet; k)\otimes H_r(Y_\bullet; \bbz) \longrightarrow CH^p(X;k+r).
\end{equation}
This pairing is functorial   both on  $X$ and on  $Y_\bullet.$
\end{proposition}
\begin{remark}
When using field coefficients for the higher Chow groups and singular homology, the slant product is compatible with the K\"unneth decomposition and cap products; see \cite{Dold-LAT}.
\end{remark}
\begin{example}
\label{ex:pairing1}
Let $X=\spec{A}$ be an  affine $\mathbb{k}$-scheme and let $B_\bullet GL_n(A)$ be the classifying space of the discrete group $GL_n(A)$, with $n>0$. The proposition above gives a slant product between the higher Chow groups of the  simplicial scheme $X \times B_{\bullet}GL_n(A)$ and the group homology of $GL_n(A)$:
\begin{equation}
    \label{eq:HCG-Hom}
   \text{-}/\text{-} \ \colon \ CH^p(X\times B_\bullet GL_n(A); k) \otimes H_r(B_\bullet GL_n(A); \bbz) \longrightarrow CH^p(X;k+r).
\end{equation}
This pairing is clearly functorial in the category of regular $\mathbb{k}$-algebras. 

When $X=\spec{A}$ is regular,  the slant product and the tautological evaluation morphism
\begin{equation}
    \label{eq:eval}
\ev \colon \mathsf{Spec}(A) \times B_\bullet GL_n(A) \longrightarrow B_\bullet GL_n,  
\end{equation}
together with a fixed class $\mbc \in CH^p(B_\bullet GL; k)$, can be used  
to define a functorial \emph{characteristic class} map
$ \gamma_\mbc \colon K_j(X) \longrightarrow CH^p(X; k+j) $
via the Hurewicz homomorphism $\mbh \colon K_j(X) = \pi_j( BGL(A)_+ ) \to H_j(B_\bullet GL(A); \bbz) $. Namely, $\gamma_\mbc$ is given by  the composition
\begin{equation}
\label{eq:gamma_c}
    \xymatrix{
K_j(X) = \pi_j( BGL(A)_+ ) \ar[r]^-{\mbh} \ar@/_1pc/[rrd]_{\gamma_\mbc} &  H_j(B_\bullet GL(A)) \ar[rd]^-{\ev^*(\mbc)/{-}}   &  \\
& &  CH^p(X;k+j).
    }
\end{equation}
\end{example}
\bigskip


\section{Constructing cycles in $B_\bullet GL_n$}
\label{sec:constr}

In order to construct special cycles in \( \sbcxeq{p}{0}{B_\bullet GL_n} \), in the sense of Definition~\ref{def:special_cycle}, that  are compatible under stabilization as \( n \to \infty\), we introduce and study a family of morphisms \( L_{n,r} \colon GL_n^{\times r} \times \Delta^r \to \mat{n}{n}\) \ that will be  used to pull-back certain algebraic cycles from \( \mat{n}{n}\). These morphisms are a key ingredient in subsequence constructions. 

\subsection{Auxiliary morphisms}
Although our primary interest lies on schemes over a field \( \Bbbk\), it is worth mentioning that the  constructions and related  discussions in this section hold over \( \bbz\). 
Define 
\(
L_{n,r} \colon \GL{n}^{\times r}\times \Delta^r \longrightarrow \mat{n}{n}
\) by 
 sending  
\(  \mbA = (A_1|\cdots | A_r) \in \GL{n}(R)^{\times r} \)
and \( \mbt=(t_0, t_1, \ldots, t_r) \in \Delta^r(R)\) to
\begin{equation}
    \label{eq:L}
L_{n,r}(\mbA, \mbt) :=  t_0 I + t_1 A_1 + t_2 (A_1A_2) +  \cdots + t_r(A_1\cdots A_r) \ \in \  \mat{n}{n}(R),
\end{equation}
where \(  I = (\delta_{ij}) \) is the identity matrix.
Given  \(  j=0,\ldots, n, \ \) let \(\{ \mbe_j \} \subset \Delta^n \) be the closed subscheme defined by the ideal 
\(\langle z_0, \ldots, z_j-1, \ldots, z_n \rangle,  \) and define \( \widehat L_{n,r} \colon GL_n^{\times r} \times \Delta^r \to \mat{n}{n} \times \Delta^r\ \)\  by\  
\( \ \widehat L_{n,r}(\mbA, \mbt) := (L_{n,r}(\mbA), \mbt).\)

\begin{remark}
\label{rem:VanEst}
One should contrast  the definition of  \( L_{n,r} \) with   parametrizations of the geodesic simplices \( \Delta(\gamma_1, \ldots, \gamma_r) \subset G/K \)  that are used in \cite{Dup-ccflat} to express the van-Est isomorphism  in terms of simplicial De Rham cohomology. We suspect that the similarity is not a mere coincidence, but underlies a De Rham realization of our geometric characteristic classes. 
\end{remark}

\begin{lemma}
\label{lem:L-smooth}
Given \( 1 \leq k \leq n-1 \),  let \( \pi_k \colon \mat{n}{n} \to \mat{n}{k}\) be the projection onto the \textbf{last} \( k \) columns. The following properties hold.
\begin{enumerate}[a.]
\item \label{it:L-a} The restriction 
\(
\widehat L_{n,r} \colon   GL_n^{\times r} \times \left( \Delta^r - \{ \mbe_0 \} \right) \ \to \ \mat{n}{n}\times \Delta^r
\)
is a smooth  morphism. 
\item \label{it:L-b}  The composition
\[
\begin{tikzcd}[column sep=small]
SL_n^{\times r} \times \left( \Delta^ r - \{ \mbe_0 \} \right) \arrow[hook]{r} &  GL_n^{\times r} \times \left( \Delta^ r - \{ \mbe_0\} \right) \arrow[r, "\widehat L"]  & 
\mat{n}{n}\times \Delta^r \arrow[r, "\pi_k" ] & \mat{n}{k}\times \Delta^r
\end{tikzcd}
\]
is smooth.
\end{enumerate}
\end{lemma}
\begin{proof}
See Appendix \ref{app:pf1}.
\end{proof}
\begin{corollary}
\label{cor:ffnr}
If  \( \scrU_{n,r} := L_{n,r}^{-1}\left( \mat{n}{n}- \{ I \} \right) 
\subset GL_n^{\times r} \times (\Delta^r - \{ \mbe_0\})\),  then the restriction
\( L_{n,r} \colon \scrU_{n,r} \to \mat{n}{n}-\{ I \} \) is faithfully flat, for all \( r \geq 1\).
\end{corollary}
\begin{proof}
Let \( I \neq M \in \mat{n}{n}(K) \) be given, where \( K \) is a field. In an algebraic closure  \(  \bar{K} \) of \( K \), choose an element \( 1\neq \lambda \in \bar{K}\) which is not an eigenvalue of \( M\). Then \( A:= \frac{1}{1-\lambda}\left( M - \lambda I \right) \) lies in \( GL_n(\bar{K}) \), while \( \mbs := (1-\lambda, \lambda )  \in \Delta^1(\bar{K}).\) 
Then \( L_{n,1}(A;\mbs) = M \), showing that \( L_{n,1}|_{\scrU_{n,1}} \) surjects onto \( \mat{n}{n}- \{ I \} \). This suffices to prove the corollary.
\end{proof}

Recall that  \(  \partial_j \colon \Delta^{r-1} \rightarrow \Delta^r  \)  and \( s_j \colon \Delta^r \to \Delta^{r-1}\) denote the co-face and co-degeneracy morphisms of \(  \Delta^\bullet \), respectively, and let \(\,  \delta_j \colon GL_n^{\times r} \to GL_n^{\times r-1}\,  \) and \(\,  \sigma_j \colon GL_n^{r-1} \to GL_n^{r}\,   \) be the face and degeneracy maps of  \(  B_\bullet GL_n. \)   The  properties below follow directly from the definitions. 
\begin{properties}
\label{rem:prop-L}
    \begin{enumerate}[i.] 
    \item Given \(  \uA \in GL_n^{\times r}(R) \) and \(  \mbs   \in \Delta^{r-1}(R)\), then
\begin{equation}
    \label{eq:rels1}
L_{n,r}(\uA; \partial_j(\mbs) ) =
\begin{cases}
A_1\cdot L_{n,r-1}(\delta_0 \uA; \mbs) &,  \text{ when } j=0\\
\ L_{n,r-1}(\delta_j \uA; \mbs) &,  \text{ when } 1 \leq j\leq r
\end{cases}
\end{equation}
\item Given \(  \uB \in GL_n^{\times (r-1)}(R) \) and \(  \mbt   \in \Delta^{r}(R)\), then
\begin{equation}
    \label{eq:rels2}
L_{n,r-1}(\uB; s_j(\mbt) ) =
  L_{n,r-1}( \sigma_j \uB; \mbt ) , \text{ when } 1 \leq j\leq r-1.
\end{equation}
\item 
Let \(  \jmath_n \colon \mat{n}{n} \hookrightarrow \mat{(n+1)}{(n+1)} \)  be the regular embedding defined by 
\begin{equation}
\label{eq:jmath}
  \begin{tikzcd}
  A \ar[r, mapsto, "\jmath_n"]  & 
  \begin{pmatrix} A & 0 \\ 0& 1 \end{pmatrix} .
\end{tikzcd}
\end{equation}
Then  the following diagram commutes.
\begin{equation}
    \label{eq:stabL}
    \xymatrix{
    GL_n^{\times r} \times \Delta^r \ar[rr]^-{L_{n,r}} \ar[d]_{\jmath_n^{\times r} \times\bone} & & \mat{n}{n} \ar[d]^{\jmath_n}\\
GL_{n+1}^{\times r} \times \Delta^r \ar[rr]^-{L_{n+1,r}} & & \mat{(n+1)}{(n+1)}
    }
\end{equation}
\end{enumerate}
\end{properties}

\subsection{Coherent families of cycles}
\label{subsec:coherent}

In this section we show how to use the properties of \( L_{n,r}\)  to construct classes in the Higher Chow Groups of \( B_\bullet GL_n \).  The first step is to characterize certain  families of algebraic cycles on the schemes \( \mat{n}{n}\) that glue together appropriately under pull-back via \( L_{n,r}\).

\begin{definition}
\label{def:coh_fam}
We say that a collection \( \{ \Gamma^p_n \mid n \geq p \} \) of algebraic cycles in \( \mat{n}{n} \)  is a \emph{coherent family of pre-cycles} of codimension \(p >0 \)  if the following conditions hold.
\begin{enumerate}[i)]
\item \label{itD:1} \( \Gamma^p_n \) is an algebraic cycle of codimension \( p \) in \( \mat{n}{n} \);
\item \label{itD:2} \( \operatorname{supp}(\Gamma^p_n) \cap \{ I \} = \emptyset \), where 
\(\operatorname{supp}(\Gamma^p_n) \)
is the support of  \(\Gamma^p_n\);
\item \label{itD:3} If \( \jmath_n \colon \mat{n}{n} \hookrightarrow \mat{(n+1)}{(n+1)} \) \, is the  embedding \eqref{eq:jmath}, then \( \jmath_n(\mat{n}{n}) \) intersects \( \Gamma^p_{n+1}\) properly and 
\[
\Gamma^p_n \ = \  \jmath_n^*(\Gamma^p_{n+1}) \equiv \jmath_{n\sharp}([\mat{n}{n}]) \centerdot \, \Gamma^p_{n+1}.  
\]
\item \label{itD:4} For each \( n \geq p \) the cycle \( \Gamma^p_n\) is invariant under the action of \( GL_n \) on \( \mat{n}{n} \) via left matrix multiplication 
\( \mu \colon GL_n\times \mat{n}{n} \to \mat{n}{n}, \ (g, B)\mapsto g\cdot B \). 
 In other words, \ 
\( g_* (\Gamma^p_n) = \Gamma^p_n\) \, for all \( g \in GL_n\).
\end{enumerate}
\end{definition}

Suppose  we are given such a  family  \( \{ \Gamma^p_n \mid \, n \geq p \} \). It follows from condition \ref{itD:2}) in the definition,  along with Corollary \ref{cor:ffnr}, that  the flat pullback 
\begin{equation}
\label{eq:flat}
 \varGamma_{n,r}^p:= L_{n,r}^* (\Gamma^p_n), \quad r\geq 1, \ \ 1\leq p \leq n, 
 \end{equation} 
 is a well-defined algebraic cycle  in \( GL_n^{\times r} \times \Delta^r \).

\begin{proposition}
\label{prop:coherent}
The cycles \( \varGamma^p_{n,r} \) satisfy the following properties. 
\begin{enumerate}[i)]
\item \( \varGamma^p_{n,r} \) has codimension \( p\) and is equidimensional and dominant over \( \Delta^r \).
\item \( (\jmath_n^{\times r}\times \bone)^* (\varGamma^p_{n+1,r} ) \ = \ \varGamma^p_{n,r}.\)
\item Given \( 0 \leq j \leq r+1 \) one has
\[
(\delta_j \times \bone)^* \varGamma_{n,r}^{p} = (\bone\times \partial_{j})^*\varGamma_{n,r+1}^p.
\]
In other words, the element \( \mathfrak{A}^p_n = (\varGamma^p_{n,r}) \ \in \ \sbcxeq{p}{0}{B_\bullet GL_n} \) is a special cycle in the sense of Definition \ref{def:special_cycle}, 
\end{enumerate}
\end{proposition}
\begin{proof}
Let \(V \) be an irreducible component of \( \operatorname{supp}(\Gamma^p_n)\). Since  \( V \cap \{ I \}= \emptyset \) (by definition), one concludes from Corollary \ref{cor:ffnr} that 
\begin{equation}
\label{eq:notempty}
\left\{ V\times (\Delta^r - \{ \mbe_0\}) \right\}\, \cap\,  \widehat{L}_{n,r}( GL_n^{\times r} \times (\Delta^t - \{ \mbe_0\} ) \ \neq \ \emptyset.
\end{equation}
One has the following commutative diagram
\[
\begin{tikzcd}
L^{-1}_{n,r}(V) =  L^{-1}_{n,r}(V) \, \cap \, \scrU_{n,r} \arrow[hook]{rr} \arrow[d] & & \scrU_{n,r} \subset GL_n^{\times r}\times (\Delta^r - \{ \mbe_0\}) \ar[d, "\widehat{L}_{n,r}"] \\
V\times (\Delta^r - \{ \mbe_0 \}) \ar[d, "pr_2"] \arrow[hook]{rr}  && (\mat{n}{n} - \{ I \}) \times (\Delta^r - \{ \mbe_0 \} ) \\
\Delta^r - \{ \mbe_0 \} & & 
\end{tikzcd}
\]
where the square is a fiber square, by definition, and the vertical arrows are flat, by Lemma \ref{lem:L-smooth}.\ref{it:L-a}. Using \eqref{eq:notempty} one concludes that every component of \( L_{n,r}^*(\Gamma^p_n) \) is equidimensional and dominant over \( \Delta^r\), thus proving statement \emph{i)}. 

The second statement follows from the definitions and \eqref{eq:stabL}, as follows.
\begin{align*}
 (\jmath_n^{\times r}\times \bone)^* (A^p_{n+1,r} ) & :=  
  (\jmath_n^{\times r}\times \bone)^*\circ  L_{n+1,r}^*(A^p_{n+1}) = 
\left\{   L_{n+1,r} \circ  (\jmath_n^{\times r}\times \bone)\right\}^*(A^p_{n+1}) \\
&\overset{\eqref{eq:stabL}}{=} 
\{\jmath_n\circ L_{n,r}\}^*(A^p_{n+1}) = L_{n,r}^*\circ \jmath_n^*(A^p_{n+1} )
= 
 L_{n,r}^* (A^p_n) = A^p_{n,r}.
 \end{align*}

For the last statement, first observe that  the identity \eqref{eq:rels1} in the case \( j = 0 \) can be expressed as the commutativity of the following diagram.
\[
\begin{tikzcd}
GL_n^{\times (r+1)} \times \Delta^r  \ar[d, "\left( \bone_{GL_n^{\times (r+1)}}\right) \times \partial_0"']  \ar[rrrr, "pr_1 \times (\delta_0 \times \bone_{\Delta^r})", " = "'] & & & & GL_n\times (GL_n^{\times r} \times \Delta^r)  \ar[d, "\bone_{GL_n} \times L_{n,r}"] \\
GL_n^{\times (r+1)} \times \Delta^{r+1} \ar[rr, "L_{n,r+1}"'] & & \mat{n}{n} & & GL_n\times \mat{n}{n} \ar[ll, "\mu"] 
\end{tikzcd}
\]
Now, it follows from property \ref{itD:4}) in Definition \ref{def:coh_fam} that \( \mu^* (A^p_n) = [GL_n] \times A^p_n \). 
Therefore, 
\begin{align}
\label{eq:partial0}
(\bone\times   \partial_0)^*(A^p_{n,r+1})  &  := (\bone \times \partial_0)^* L^*_{n,r+1}(A^p_n) =
(\bone_{GL_n} \times L_{n,r})^* \mu^* (A^p_n)  \\
& =
(\bone_{GL_n} \times L_{n,r})^* ([GL_n]\times A^p_n) = [GL_n]\times L_{n,r}^*(A^p_n) 
\notag \\
&  =
[GL_n]\times A^p_{n,r}. \notag
\end{align}
On the other hand, since \( \delta_0 \colon GL_n^{\times (r+1)} \to GL_n^{\times r}\) is the projection onto the last \( r \) factors, one concludes that under 
\( \delta_0 \times \bone \colon GL_n^{\times (r+1)} \times \Delta^r \to GL_n^{\times r}\times \Delta^r \) one has 
\( (\delta_0 \times \bone)^* (A^p_{n,r}) = [GL_n]\times A^p_{n,r},\)
and hence 
\( (\delta_0 \times \bone)^\flat (A^p_{n,r}) = (\bone \times \partial_0)^*(A^p_{n,r+1}),\)
according to \eqref{eq:partial0}.

Now, when \(1\leq j \leq r+1\), property \eqref{eq:rels1} simply states that 
\( L_{n,r+1} \circ (\bone\times \partial_j)  =  L_{n,r}\circ (\delta_j \times \bone) \). Therefore, one obtains
\( (\bone\times \partial_j) ^*(A^p_{n,r+1}) = (\delta_j \times \bone)^*(A^p_{n,r}) \) and this concludes the proof. 

\end{proof}


%

\section{Coherent families of determinantal schemes}
\label{sec:det_fam}

This section introduces the determinantal schemes that are used in the construction of the explicit Chern cycles in \S \ref{sec:Ccycles}. We study their ideals of definition and highlight  key properties and relations amongst them, following closely the notation in \cite{GortzWed-AG}. The reader is encouraged to consult \cite{BrunsHerz-CMrings}, \cite{Ful-IT} and \cite{Ful-det} for further results and details. 


\subsection{Preliminaries on determinantal schemes}

Let   \(  \scrE \) and \(  \scrF \) be locally free $\scrO_S$-modules of finite ranks \(  m \) and \(  n \), respectively, over a base-scheme $S.$ Given $r \geq 0$, let \( \detsch{r}{\scrE}{\scrF}\  \in  \osch{S} \) \ be  the determinantal scheme representing the functor \(  \fdetsch{r}{\scrE}{\scrF}  \colon \osch{S} \to \sets\) that sends \(  f \colon T \to S \) to 
\[  
\detsch{r}{\scrE}{\scrF} (T) := \{ u \in  \Hom{\scrO_S}{f^*\scrE}{f^*\scrF} \mid  \wedge^{r+1} u  = \mbzero \}.
\]
These schemes satisfy the following properties; see \cite[\S 16]{GortzWed-AG}.

\begin{enumerate}[i)]
   \item \(  \detsch{r}{\scrE}{\scrF} \) is a closed subscheme of \(  \inthom{\scrE}{\scrF} \) and is affine over \(  S; \)
   \item \(  \detsch{r}{\scrE}{\scrF} \) is separated and of finite presentation over \(  S \);
   \item If \(  0\leq r \leq \text{min}\{ m, n\} \) then \(  \detsch{r}{\scrE}{\scrF} \) is flat of relative dimension \(  r(m+n-r)  \) over \(  S \);    \item Assume that \(  S \) is locally noetherian. If \(  S \)  is \emph{reduced} (resp. \emph{normal},  \emph{irreducible} or \emph{Cohen-Macaulay}), then so is \(  \detsch{r}{\scrE}{\scrF} \).
\end{enumerate}

 Let \( \mathbf{x} = (x_{ij}) \)  denote an \(n\times k\) matrix of variables and let \(  \mat{n}{k}\cong\ \bba^{nk} \) be the affine space 
\( \spec{\bbz[\mathbf{x}]} \)  of  \(  n\times k \) matrices. As usual,   the general linear group \(   \GL{n} \) is realized as  an open subscheme of \( \mat{n}{n}\).   
\begin{example}
\label{exmp:detsch_1}
Consider \(  S= \spec{\bbz},\  \scrE= \scrO^{\oplus k}_{\spec{\bbz}} \)\ and\ \(  \scrF= \scrO^{\oplus n}_{\spec{\bbz}} \). Hence, one can identify \(    \inthom{\scrO^{\oplus k}_{\spec{\bbz}}}{ \scrO^{\oplus n}_{\spec{\bbz}} } \) with \(\mat{n}{k} \), and for \(0 \leq r \leq \text{min}\{ k, n \} \) the scheme \(  \detsch{r}{\scrE }{ \scrF } \) is the usual determinantal subscheme \(  D^r_{n,k} \subset \mat{n}{k}\) whose ideal is generated by the \(  (r+1)\times (r+1) \) minors of the \(  n\times k \) matrix of indeterminates \( \mathbf{x}\). See \cite{BrunsHerz-CMrings}.
\end{example}
\begin{example}
\label{exmp:detsch_2}
 Let  
\(  \scrQ_{n-1} \)  be the  dual of the universal quotient bundle of rank $n-1$  over \( \bbp^{n-1} \), fitting in  the   exact sequence 
\[
0 \to  \scrQ_{n-1} \to  \scrO_{\bbp^{n-1}}^{\oplus n} \to \scrO_{\bbp^{n-1}}(1) \to 0.
\]
Then, 
\(  \inthom{ \scrO_{\bbp^{n-1}}^{\oplus k}}{\scrQ_{n-1}} \) is a closed subscheme of \(  \bbp^{n-1} \times \mat{n}{k} \) and,   for all \(  k\geq  1 \leq r \leq  \text{min}\{ k, n-1\} \)  the determinantal scheme \(  \detsch{r}{\scrO_{\bbp^{n-1}}^{\oplus k}}{\scrQ_{n-1} } \) is an  integral, normal, Cohen-Macaulay closed subscheme of \(  \bbp^{n-1}\times \mat{n}{k} \),   flat of relative dimension \(  r(k+n-1-r) \) over \(  \bbp^{n-1} \).  
\end{example}

\begin{definition}
\label{def:main_det}
Let \(  \pi_k \colon \mat{n}{n} \to \mat{n}{k} \) denote the  projection onto the \textbf{last} \(  k \)  columns of an \(  n\times n \) matrix, and consider $1\leq p, q, j \leq n$.   
\begin{enumerate}[i.]
\item \label{it:bbd} Define
\begin{equation}
\label{eq:bbd}
\bbd^p_n := \pi_{n-p+1}^{-1}( D^{n-p}_{n, \, n-p+1}) \subset \mat{n}{n},
\end{equation}
where $D^{n-p}_{n,\, n-p+1} \subset \mat{n}{n-p+1}\) is the subscheme introduced in Example \ref{exmp:detsch_1}.
It follows that \(\bbd^p_n\)  is an irreducible, normal and Cohen-Macaulay closed subscheme of codimension \(  p\).
\item Define
\begin{equation}
\label{eq:theta}
\theta_n^{q} := (1\times \pi_q)^{-1}(\qbdlp{q}{n-1} )
\cong \mat{n}{(n-q)}\times   \qbdlp{q}{n-1}\  \subset\  \bbp^{n-1}\times \mat{n}{n}.
\end{equation}
This is a regular, closed subscheme of codimension \(  q. \)
\item Finally, define 
\begin{equation}
\label{eq:S}
\bbs_{n,j} := (1\times \pi_{j})^{-1} \left(    \detsch{j - 1}{\scrO_{\bbp^{n-1}}^{\oplus j}}{\scrQ_{n-1}}   \right) \ \subset \ \bbp^{n-1}\times \mat{n}{n}.
\end{equation}
This is a closed subscheme of codimension \(  n \)   that satisfies all properties listed in Example \ref{exmp:detsch_2}. It is convenient to set
\( \bbs_{n,0} = \bbs_{n,n+1} = \emptyset\).

\end{enumerate}
\end{definition}

\subsection{Ideals and intersections}
Our goal is to exhibit simple presentations for the ideals defining  \( \bbd^p_n\, , \, \theta^q_n\, , \, \bbs_{n,j} \) and   describe   relations between them. For brevity,  set $\mbu = (u_1, \ldots, u_n)$ and let $\mbx^k =( x_{1,k}, \ldots, x_{n,k} )^T$ be the $k$-th column of the  matrix of variables $\mbx = (x_{ij})$.

\begin{definition}
\label{def:ideals}
  Given \( I = (1\leq i_p < i_{p+1} < \cdots < i_n \leq n ) \in \Lwedge^{n-p+1}[n] \) and \( k\in [n]\),  denote 
\begin{align*}
\mbm_{p, I} & := \det \left(x_{i_r, \, s} \right)_{p \leq r, s \leq n } \, \in\, \bbz[\mbx],
\quad   \text{and}  \\
\mbu\cdot \mbx^k & := u_1 x_{1,k} + \cdots + u_n x_{n,k} \, \in \, \bbz[\mbx,\mbu].
\end{align*}
Define the following ideals
\begin{align}
\fra_p & := \la \mbm_{p, I} \mid I \in \Lwedge^{n-p+1}[n] \ra\ \subset\  \bbz[\mbx] \\
\frb_p & := \la \mbu\cdot \mbx^{p+1}, \ldots, \mbu \cdot \mbx^n \ra\   \subset\  \bbz[\mbu,\mbx] \\
\frA_p &:=  \bbz[\mbu,\mbx]\cdot \fra_p + \frb_p \  \subset\  \bbz[\mbu,\mbx]\\
\varSigma_p &:= \bbz[\mbu,\mbx]\cdot \fra_p + \frb_{p-1}\ \subset \ \bbz[\mbu,\mbx].
\end{align}
Note that \( \fra_p \subset \fra_{p+1}\) and \( \frb_{p+1} \subset \frb_{p}\), when \(1\leq p <n \). 
\end{definition}

To simplify notation in what follows, we write \( A= \bbz[\mbx] \) so that \( \bbz[\mbu, \mbx]\) is the polynomial ring \(A[\mbu] \) over \( A\),  and 
\( \bbp^{n-1}_A = \proj{A[\mbu]} = \bbp^{n-1}\times \mat{n}{n}\).

\begin{lemma}
\label{lem:tech1}
Given \( 1 \leq p \leq n\), the following holds.
\begin{enumerate}[\rm i)]
     \item \label{it:tech1} The ideal \( \fra_p \subset A\) defines the subscheme  
     \begin{equation}
     \label{eq:ideal_D}
      \bbd^p_n= \spec{A/\fra_p} \subset \mat{n}{n}.
      \end{equation}
     \item \label{it:tech2} The homogeneous ideal \( \frb_q \subset A[\mbu]\) defines the subscheme 
      \begin{equation}
      \label{eq:ideal_theta}  
      \theta^{n-q}_n   =  \proj{A[\mbu]/\frb_{q}}    \subset \bbp^{n-1}\times \mat{n}{n}.
      \end{equation}
    \item \label{it:tech3}  The homogeneous ideal $\varSigma_j \subset A[\mbu] \) is prime and defines the subscheme
     \begin{equation}
     \label{eq:ideal_S} \bbs_{n,j} = \proj{A[\mbu]/\varSigma_j} \subset \bbp^{n-1}\times \mat{n}{n}.  
     \end{equation}
\end{enumerate}
\end{lemma}
\begin{proof}
All  statements should be evident from the definition of the subschemes; see  \cite{BrunsHerz-CMrings} and  \cite[\S 16]{GortzWed-AG}.
\end{proof}

\begin{remark}
\label{rem:pts-bbp}
A matrix 
\(  A\in \mat{n}{n}(R) \) lies in \(  \bbd^p_n(R) \) iff \(  A^p\wedge A^{p+1}\wedge \cdots \wedge A^n = 0 \in \Lwedge^{n-p+1}(R^{\, \oplus n}), \) where \( A^j \) denotes the \(j\)-th column of  \( A\). 
\end{remark}

The following observation also follows directly from the definitions.

\begin{lemma}
\label{lem:stab_D}
Let  \( \imath_{n-1}\colon \bbp^{n-1} \hookrightarrow \bbp^n\) be the regular embedding  
\begin{equation*}
\label{eq:stab}
  \begin{tikzcd}
\left[u_1:\, \cdots\, : u_n  \right]   \ar[r, mapsto, "\imath_{n-1}"] & 
\left[u_1:\, \cdots\, : u_n:0 \right],
\end{tikzcd}
\end{equation*}
and \(\jmath_n \colon \mat{n}{n}\hookrightarrow \mat{(n+1)}{(n+1)} \) be as in \eqref{eq:jmath}.
Then,   for  \(  1\leq p \leq n \), one has:
\[
\jmath_n^{-1}( \bbd^p_{n+1}) = \bbd^p_n, \quad (\imath_{n-1} \times \jmath_n)^{-1}\left( \theta^p_{n+1}  \right) = \theta^p_n , \quad \text{and} \quad  (\imath_{n-1} \times \jmath_n)^{-1}\left( \bbs_{n+1,r}\right) = \bbs_{n,r} .
\]

\end{lemma}

\begin{proposition}
\label{prop:intersect}
 If $1\leq p \leq n-1$, then
\begin{equation}
\label{eq:inter_ideals}
 \mathfrak{A}_p = \varSigma_{p} \cap \varSigma_{p+1}.
 \end{equation}
 Furthermore, if  $\mathfrak{m}_p$  is  the maximal ideal in the local ring $A[\mbu]_{\varSigma_{p}}$ then
   \begin{align}
     \mathfrak{m}_p & =  \fra_p\cdot A[\mbu]_{\varSigma_{p}} + \frb_{p} \cdot A[\mbu]_{\varSigma_{p}} \quad \text{and} \quad \label{eq:id1} \\
    \mathfrak{m}_{p+1} & = \fra_p \cdot A[\mbu]_{\varSigma_{p+1}} +\frb_{p} \cdot A[\mbu]_{\varSigma_{p+1}}. \label{eq:id2}
  \end{align}
\end{proposition}

\begin{proof}
For simplicity,  write \( \hat\fra_p := A[\mbu]\cdot \fra_p \subset A[u]\), so that \(  \frA_p = \hat \fra_p + \frb_p\), 
\( \varSigma_p = \hat \fra_p + \frb_{p-1}\) and \( \varSigma_{p+1} = \hat \fra_{p+1} + \frb_{p}\).
Since  \( \fra_p \subset \fra_{p+1}\) and \( \frb_{p+1} \subset \frb_{p}\), one concludes that \( \, \frA_p \subset \varSigma_{p} \cap \varSigma_{p+1} \).

On the other hand, 
\begin{align}
\label{eq:inter2}
 \varSigma_p \cap \varSigma_{p+1} & = (  \hat \fra_p + \frb_{p-1} ) \cap (\hat \fra_{p+1} + \frb_{p} )  \notag \\
 & = \hat \fra_p \cap ( \hat \fra_{p+1} + \frb_p ) + \frb_{p-1}\cap \frb_p + \frb_{p-1}\cap \hat \fra_{p+1}\\
& \subset (\fra_p+ \frb_p ) + \frb_{p-1}\cap \hat \fra_{p+1} \ = \ 
\frA_p \ + \ (\hat\fra_{p+1} \cap \frb_{p-1}).
\notag
\end{align}
Hence,   identity \eqref{eq:inter_ideals} is a consequence of the following claims. 
\begin{enumerate}[{\sc Claim} 1.]
\item \( \hat \fra_{p+1} \cdot \frb_{p-1} \subset \frA_p \)
\item  \( \hat \fra_{p+1} \cap \frb_{p-1} =  \hat \fra_{p+1} \cdot \frb_{p-1} \).
\end{enumerate}

To prove Claim 1, first observe that \( \frb_{p-1} = \la \mbu\cdot \mbx^p \ra + \frb_p  \). It follows that  \( \frb_{p-1} \cdot \hat \fra_{p+1} \subset  \la \mbu\cdot \mbx^p \ra \cdot \hat \fra_{p+1} + \frA_p\) and we just need to show that 
\begin{equation}
\label{eq:tricky}
(\mbu\cdot \mbx^p)\, \mbm_{p+1, I}\, \in\, \frA_p, \quad \text{for all} \quad I \in \Lwedge^{n-p}[n].
\end{equation}
Let us start with the case \( I_0 = ( p+1< \cdots  <n ) \). Consider the \( (n+1)\times (n+1) \) matrix
\begin{equation}
\label{eq:detM}
M_{p, I_0} =
\begin{pmatrix}
1         & \cdots & 0         & x_{1,p+1}& \cdots &x_{1,n} &x_{1,p} & \\ 
\vdots &            & \vdots & \vdots       &            &\vdots    & \vdots  & \\
0         & \cdots & 1         & x_{p,p+1}& \cdots &x_{p,n} &x_{p,p} & \\
0         & \cdots & 0         & x_{p+1,p+1}& \cdots &x_{p+1,n} &x_{1,p} & \\ 
\vdots &            & \vdots & \vdots       &            &\vdots    & \vdots  & \\
0         & \cdots & 0         & x_{n,p+1}& \cdots &x_{n,n} &x_{n,p} & \\
u_1     & \cdots & u_p     & \mbu\cdot \mbx^{p+1}  & \cdots & \mbu\cdot \mbx^n  &\mbu\cdot \mbx^p & \\
\end{pmatrix}.
\end{equation}
Expanding the determinant along the bottom row gives
\begin{equation}
\label{eq:detM2}
\det(M_{p,I_0}) = (-1)^{p+1} \underbrace{\left( \sum_{\ell=1}^p  u_\ell \, \mbm_{p, \{ \ell\} \cup I_0} \right)}_{\in\  \hat\fra_p} \ + \ 
\underbrace{\left( \sum_{\ell=p+1}^n (\mbu\cdot \mbx^\ell)\, F_\ell \right)}_{\in \ \frb_{p} }\ + \ (\mbu\cdot \mbx^p) \, \mbm_{p+1,I_0},
\end{equation}
where \( F_\ell \) is a polynomial in \( A = \bbz[\mbx]\).  Observe that the bottom row is a linear combination with coefficients in \( \bbz[\mbu]\)  of the other rows, and hence \( \det(M_{p,I_0}) = 0\). It follows that 
\( (\mbu\cdot \mbx^p) \, \mbm_{p+1,I_0} \in \hat \fra_{p} + \frb_{p}= \frA_p \).

Given an arbitrary $I \in \wedge^{n-p+1}[n]$, one can construct a similar matrix $M_{p,I}$ from which an identity corresponding to \eqref{eq:detM2} follows when $I_0 $ is  replaced by $I$. This shows that $(\mbu\cdot \mbx^{p})\cdot \hat\fra_{p+1} \subset \frA_p $
and proves {\sc Claim~1}, according to \eqref{eq:tricky}.

Now, observe that $(\mbu\cdot \mbx^p, \ldots, \mbu\cdot \mbx^{n})$ is a regular sequence in $A[\mbu]/\hat \fra_{p+1}$. Since these are the generators of $\frb_{p-1}$, it follows that $\mathsf{Tor}^1_{A[\mbu]}(A[\mbu]/\frb_{p-1}, A[\mbu]/\hat\fra_{p+1}) = 0$ and hence   {\sc Claim 2} follows. This concludes the proof of the first assertion in the proposition.

We now prove \eqref{eq:id1} and \eqref{eq:id2}. By definition, the maximal ideals are given by 
$ \mathfrak{m}_p = \fra_p\cdot A[\mbu]_{\varSigma_{p}} + \frb_{p-1} \cdot A[\mbu]_{\varSigma_{p}}$
and 
$\mathfrak{m}_{p+1} = \fra_{p+1}\cdot A[\mbu]_{\varSigma_{p+1}} + \frb_{p} \cdot A[\mbu]_{\varSigma_{p+1}}.$
Therefore, it suffices to show that
\begin{equation}
\label{eq:IJ1}
\frb_{p-1} \cdot A[\mbu]_{\varSigma_{p}} \subset 
(\hat\fra_{p} + \frb_{p})\cdot A[\mbu]_{\varSigma_{p}} = \frA_p\cdot A[\mbu]_{\varSigma_{p}}
\end{equation}
and  
\begin{equation}
\label{eq:IJ2}
\hat\fra_{p+1} \cdot A[\mbu]_{\varSigma_{p+1}} \subset 
(\hat\fra_{p} + \frb_{p} )\cdot A[\mbu]_{\varSigma_{p+1}} = \frA_p \cdot   A[\mbu]_{\varSigma_{p+1}}.
\end{equation}

Note that  \( \mbm_{p+1, I_0} \notin \varSigma_p\), and hence \( \mbm_{p+1,I_0} \) becomes a unit in the local ring \( A[\mbu]_{\varSigma_{p}} \). On the other hand,  \eqref{eq:tricky} shows that that \( (\mbu\cdot \mbx^p)\,  \mbm_{p+1,I_0} \in \frA_p\). Therefore \( \mbu \cdot \mbx^p \) lies in \(\frA_p \cdot A[\mbu]_{\varSigma_{p}}
 \), thus proving \eqref{eq:IJ1}.

Similarly,   $(\mbu\cdot \mbx^{p}) \notin \varSigma_{p+1}$ and  
$\mbu\cdot \mbx^{p+1}$ is a  unit $A[\mbu]_{\varSigma_{p+1}}.\, $  Once again,  invoking \eqref{eq:tricky}  one concludes that  
\(\mbm_{p+1,I} \in \frA_p \cdot A[\mbu]_{\varSigma_{p+1}}\) for all \( I \in \Lwedge^{n-p}[n] \), and \eqref{eq:IJ2} follows.
\end{proof}

 \begin{corollary}
 \label{cor:sch_int}
  Let \(  \rho \colon \bbp^{n-1}\times \mat{n}{n} \to \mat{n}{n} \) denote the projection.  Given $0 < p < n$, the closed subschemes 
  \(  \theta^{n-p}_n \) and \( \rho^{-1}(\bbd^{p}_n )\) 
of \( \bbp^{n-1}\times \mat{n}{n} \)
 intersect properly, and  \(\theta^{n-p}_{n}\, \cap\, \rho^{-1}(\bbd^{p}_n)  = \bbs_{n, p} \ \cup \ \bbs_{n, p+1}.\) Furthermore, over a field \( \Bbbk \) both irreducible components have geometric multiplicity \(  1 \), and one has the following identity of algebraic cycles:
\begin{equation}
\label{eq:inter}
[\theta^{n-p}_{n | \Bbbk}] \centerdot  \rho^*[\bbd^{p}_{n |  \Bbbk}] \ =  [\bbs_{n,p  | \Bbbk}] + [ \bbs_{n,p+1 | \Bbbk}].
\end{equation}

 \end{corollary}

\begin{proof}
It follows from Lemma \ref{lem:tech1}.\ref{it:tech1} and \ref{lem:tech1}.\ref{it:tech2} that 
\begin{align*}
\theta^{n-p}_n \cap \rho^{-1}(\bbd^{p}_n ) & =
\proj{A[\mbu]/\frb_{p}} \cap\proj{A[\mbu]/A[\mbu]\cdot \fra_{p}}    \\
&  = 
\proj{ A[\mbu]/ \{ \frb_p + A[\mbu]\cdot \fra_{p} \} } =:
\proj{ A[\mbu]/ \mathfrak{A}_p }.
\end{align*}
On the other hand, \eqref{eq:inter_ideals}  shows that 
\[  
\proj{A[\mbu]/\mathfrak{A}_p } = \proj{A[\mbu]/\varSigma_p \cap \varSigma_{p+1}} = \bbs_{n,p}\cup \bbs_{n,p+1},
\]
and the  statement on multiplicities comes from  \eqref{eq:id1}, \eqref{eq:id2} and  \cite[Prop.~8.2]{Ful-IT}.
\end{proof}

\section{The Chern cycles}
\label{sec:Ccycles}

In this section we  show that, after base change to a field $\mathbb{k}$,  the determinantal schemes introduced in the previous section assemble to give explicit cycles
\( \mathfrak{C}^p_n\in \sbcxeq{p}{0}{B_\bullet GL_{n, \mathbb{k}}} \)   that represent the $p$-th Chern class \(  \mbc_p \in CH^p(B_\bullet GL_{n, \mathbb{k}};0)\) of the universal \(  n \)-plane  bundle, for $1 \leq p \leq n$.
These cycles are  compatible with the inclusions 
$GL_n \hookrightarrow GL_{n+1}$ and hence they represent the generators of the Chow ring of $B_\bullet GL_{\mathbb{k}}.$

Let us summarize the relevant properties obtained thus far.

\begin{proposition}
The subschemes \( \bbd^p_n \subset \mat{n}{n}\) satisfy the following properties.
\begin{enumerate}[i.]
\item \( \operatorname{codim}(\bbd^p_n) = p \);
\item \( \jmath_n^{-1}( \bbd^p_{n+1}) = \bbd^p_n\);
\item \( \bbd^p_n \cap \{ I \}  \ = \ \emptyset\);
\item \( \bbd^p_n\) is invariant under the action of \( GL_n\) on \( \mat{n}{n} \) via left matrix  multiplication.
\end{enumerate}
\end{proposition}
\begin{proof}
See Definition \ref{def:main_det} for the first statement and Lemma \ref{lem:stab_D} for the second one. The third and fourth statements follow directly from Remark \ref{rem:pts-bbp} and definitions. 
\end{proof}

\begin{corollary}
\label{cor:coh_famD}
For fixed \(p>0 \), the collection \( \{ [\bbd^p_n] \mid n\geq p \} \) forms  a coherent family of pre-cycles of codimension \( p \) in the sense of Definition \ref{def:coh_fam}.
\end{corollary}

This corollary suffices to apply Propositon \ref{prop:coherent} to  the family \( \{ [\bbd^p_n ]\, | \, n\geq p\} \). However, in this particular case one  can derive stronger properties directly at the level of subschemes, as explained next. This may be useful to study further realizations of Chern classes, from higher algebraic \( K\)-theory into other cohomology theories, since the underlying subschemes used have  many amenable properties.

\begin{definition}
\label{def:}
Given \(  r\geq 1  \) and \(  1\leq p \leq n \), define closed subschemes
\begin{equation}
\label{eq:Cs}
C^{p}_{n,r} := L_{n,r}^{-1}(\bbd^p_n)\ \subset\ GL_n^{\times r} \times \Delta^r,
\end{equation}
and note that \( [C^p_{n,r}] = L_{n,r}^*[\bbd^p_n] \).
\end{definition}
\begin{proposition}
\label{prop:cycles2}
The subschemes \( C^p_{n,r}\) satisfy the following properties. 
\begin{enumerate}[i.]
   \item The functor of points \(  R\mapsto C^{p}_{n,r}(R) \)  is given by  
   \begin{equation}
\label{eq:C(R)}
\hche{p}{r}{n}(R) := \{ (\mbA, \mbt) \mid
L^p_{n,r}(\mbA, \mbt)\wedge \cdots \wedge L^{n}_{n,r}(\mbA, \mbt) = \mbzero \},
\end{equation}
where \(L^j_{n,r}(\mbA, \mbt)\) denotes the \(  j \)-the column of the matrix \(L_{n,r}(\mbA, \mbt) \).

   \item \label{it:2a} \(  \hche{p}{r}{n} \) is Cohen-Macaulay, normal, irreducible and has codimension \(  p; \)
   \item \label{it:2b} \(  \hche{p}{r}{n}  \) is equidimensional and dominant over \(  \Delta^r \).
   \item \label{it:2c} For $1\leq p \leq n$ one has $(\jmath_n^{\times r} \times\bone)^{-1}\left(C^{p}_{n+1, r}\right) = C^{p}_{n,r};$ see \eqref{eq:stabL}.
\end{enumerate}
\end{proposition}
\begin{proof}
The first assertion follows directly from the definitions of \(  L_{n,r}\)  and \(\bbd^p_{n}\).

The notation used below comes from  Figure \ref{fig:compo} (Appendix \ref{app:pf1}). As a consequence of  Lemma \ref{lem:L-smooth} one sees that 
 the restriction 
\[  \widehat L_{n,r}  \colon  (F_{n,r}\times p_o)^{-1}\left( \jmath^{-1} (\scrU_{n,r})  \right)
\longrightarrow \mat{n}{n}- \{ I \}
\]
is also smooth.
On the other hand, the inclusion  \(  \bbd^q_n \subset \mat{n}{n}- \{ I \} \)  implies that 
\(  C_{n,r}^{p} = ( \widehat L_{n,r})^{-1}(\bbd^p_n) \), and  since \(  \widehat L_{n,r} \) is smooth,  it follows from   Example \ref{exmp:detsch_1} that  \(  C^{p}_{n,r} \) satisfies all properties in assertion \ref{it:2a}.

The third assertion can be shown using the arguments in Proposition \ref{prop:coherent}. Alternatively, we can use Lemma \ref{lem:smooth} where we show that \(  \mu_{n,r|_{\scrU_{n,r}}} = pr_1\circ \hatmu_{n,r} \), where \(  \hatmu_{n,r} \) is  smooth and \(  pr_2 \circ \hatmu_{n,r} \) is also smooth and dominant. As a result, one concludes that \(C^{p}_{n,r}\) is equidimensional and dominant over \(  \Delta^r. \)

Given $1\leq k \leq n$ and $A \in \mat{n}{n}(R)$, the  $k$-th column $\jmath_n(A)^k$ of $\jmath_n(A)$ is precisely the augmented column $\begin{pmatrix} A^k \\ 0 \end{pmatrix}$, where \(A^k\) is  $k$-th column of $A$, while \( \jmath_n(A)^{n+1} = \mbe_{n+1}\). In particular,
$\jmath_n(A)^p \wedge \cdots\wedge \jmath_n(A)^{n}\wedge \jmath_n(A)^{n+1} =0  $ if and only if 
$A^p \wedge A^{p+1}\wedge \cdots\wedge  A^n =0  $.
From this identity one concludes that $\jmath_n^{-1}(\bbd_{n+1}^k) = \bbd^k_n,$ for $1\leq k \leq n$.
 Assertion 
\ref{it:2c} now follows from this fact and definitions.

\end{proof}
%

%

%
\begin{definition}
\label{def:Chern_cyc}

Fix  \(  1 \leq p \leq n \),  and define
\begin{equation}
\label{eq:CCp}
\mathfrak{C}^p_n := ([C^{p}_{n,r}] ) \in \sbcxeq{p}{0}{B_\bullet GL_n} = \prod_r \sbcxeq{p}{r}{\glt{n}{r}}. 
\end{equation}
The element \( \mathfrak{C}^p_n \) is a \emph{special cycle}  in the sense of Definition~\ref{def:special_cycle}, according to Corollary \ref{cor:coh_famD} or Proposition \ref{prop:cycles2},    and \S \ref{subsec:coherent}. Hence it represents a class in \( \shcg{p}{0}{B_\bullet GL_n} \).\end{definition}

The main result  asserts that the \( \mathfrak{C}^p_n\) are the desired Chern cycles.

\begin{theorem}
\label{thm:main-T}
For each $1\leq p\leq n$, the element \( \mathfrak{C}^p_n \in \sbcxeq{p}{0}{\bgl{n}} \) is a cycle in the higher Chow complex of $\bgl{n}$ that represents the \(  p \)-th  Chern class \(  \mbc_p \in \shcg{p}{0}{\bgl{n}}  \) of the universal \(  n \)-plane bundle \(  \bbe_{n\, \bullet} \)  over \( \bgl{n}. \) Furthermore, the cycles \( \mathfrak{C}^p_n  \)  are compatible under  pull-back via the inclusions \(  \jmath_n \colon GL_n \hookrightarrow GL_{n+1} \).
\end{theorem}

\begin{remark}
\label{rem:FCC}
The arguments used in Example \ref{ex:pairing1}, along with this theorem, directly provide  functorial Chern classes
\[
c_{p,j} \colon K_j(X)\to CH^p(X;j)
\]
for \( X\) affine and smooth over \( \Bbbk\). Therefore, using homotopy invariance one can apply the Jouanolou-Thomason trick to define Chern classes for arbitrary smooth schemes over \( \Bbbk\).
\end{remark}

\subsection{Proof of Theorem \ref{thm:main-T}}

We first  introduce a  family of special cycles that   represent  powers of the first Chern class of  \( \scrO_{\bbp(\bbe_n)}(1) \)  in the higher Chow ring of \( \bbp_\bullet(\bbe_n) \), the projectivization of the universal \( n\)-plane bundle \( \bbe_n \to B_\bullet GL_n\); see Example \ref{exmp:bar}.  The main ingredients are the subschemes \( \theta^q_n \subset \bbp^{n-1}\times \mat{n}{n} \) introduced in Definition \ref{def:main_det}.

\begin{remark}
\label{rem:delta}
For clarity, denote the face and degeneracy maps of \( \bbp_\bullet(\bbe_n)\)  by  
\[ 
\hat \delta_j \colon \bbp^{n-1}\times GL_n^{\times r} \to \bbp^{n-1}\times  GL_n^{\times r-1} \quad \text{and}\quad  \hat\sigma_j \colon \bbp^{n-1}\times  GL_n^{r-1} \to \bbp^{n-1}\times  GL_n^{r},
\]
respectively.
\end{remark}
\begin{definition}
\label{def:theta}
Given \(  r\geq 1  \) and \(  1\leq q \leq n \), define 
\begin{equation*}
\theta^{q}_{n,r} := (\bone\times L_{n,r})^{-1} \left( \theta^{q}_n \right)\ \subset\  \bbp^{n-1} \times GL_n^{\times r} \times \Delta^r. 
\end{equation*}
where \( \theta^q_r\) is defined in \eqref{eq:theta}.
\end{definition}

\begin{proposition}
\label{prop:cycles3}
The subschemes \( \theta^q_{n,r}\) satisfy the following properties.
\begin{enumerate}[1.]
   \item The functor of points  \(  R\mapsto \theta^{q}_{n,r}(R) \) is given by
   \begin{equation}
\label{eq:T(R)}
\theta^{q}_{n,r}(R):= \{ ([u], \mbA, \mbt) \mid u \cdot L^j(\mbA, \mbt) = 0, \ j=n-q+1,\ldots, n \},
\end{equation}
where \(  [u] \) denotes the class of \(  u \in (\bba^n - \{ 0 \})(R) \) in \(  \bbp^{n-1}(R) \), and \(  u\cdot L^j(\mbA, \mbt)  \) denotes matrix multiplication when we see \(  u \) and \(  L^j(\mbA, \mbt)  \) as  elements in \(  \mat{1}{n}(R) \) and \(  \mat{n}{1}(R), \) respectively.
   \item \label{it:3a} \(  \theta^{q}_{n,r} \) is an integral subscheme of codimension \(  q; \)
   \item \label{it:3b} \( \theta^{q}_{n,r} \) is smooth over \( \bbp^{n-1}\times \Delta^r \), and equidimensional and dominant over \(  \Delta^r \).
   \item For $1\leq q \leq n$ one has $(\bone\times \jmath^{\times r}_n\times\bone)^{-1}\left(\theta^{q}_{n+1,r}\right) = \theta^{q}_{n,r};$ see \eqref{eq:stabL}.
\end{enumerate}
\end{proposition}
\begin{proof}
The description of the functor of points is a rephrasing of the definition. 

The other properties   follow by proving their counterparts for the cone \(  C\theta^{q}_{n,r}  \subset (\bba^n -\{ 0 \}) \times GL_n^{\times r} \times \Delta^r \) 
over \( \bbp^{n-1} \times GL_n^{\times r} \times \Delta^r.    \)
This is a direct consequence of  Lemma~\ref{lem:smooth}.\ref{it:Lb}, Figure \ref{fig:compo} and definitions. 
\end{proof}

Finally, for  fixed  \(  1 \leq q \leq n \),  denote
 \begin{equation}
\label{eq:power}
\Theta^q_n := ( [ \theta^{q}_{n,r}]) \in \sbcxeq{q}{0}{\bbp_\bullet(\bbe_n)}. 
\end{equation}

\begin{lemma} 
\label{lem:aux2}
Given $1\leq q \leq n $ and $0\leq j \leq r+1 $, the following identity holds.
\begin{align}
   \label{eq:id_theta}  (\delta_j \times \bone)^* [\theta^q_{n,r}] & = (\bone \times \partial_{j})^*[\theta_{n,r+1}^q].
\end{align}
In other words, $\Theta^q_n$ is a \emph{special cycle} in the sense of Definition \ref{def:special_cycle}.
\end{lemma}
\begin{proof}
Consider the morphisms
\[
\begin{tikzcd}[column sep=1ex]
&  \arrow[dl, "\hat \delta_j\times\bone"']  (  \bbp^{n-1}\times GL_n^{\times r+1})\times \Delta^r  \arrow[dr, " \bone
\times \partial_j"]  & \\
(\bbp^{n-1}\times  GL_n^{\times r})\times \Delta^r   &    &  (\bbp^{n-1}\times GL_n^{r+1})\times \Delta^{r+1} .
\end{tikzcd}
\]
%
%
Once again, using \eqref{eq:rels1} and \eqref{eq:rels2} we only need to consider consider the case $j=0$ and   associated functor of points.
 %
 %
 By definition,  a point
 \[
 ([\mbu], \uA, \mbs) \in \left( \bbp^{n-1}\times GL_n^{\times r+1}\times \Delta^r \right)(R)
 \] 
 lies in 
 $(\bone\times \partial_0)^{-1}(\theta^{q}_{n,r+1})$
 if and only if
 $\mbu\cdot L^k(\uA, \partial_0(\mbs))  = 0,$ for all \( k = n-q+1, \ldots, n\).
 From \eqref{eq:rels1} one concludes that this occurs whenever
 \begin{equation}
 \label{eq:last0}
 (\mbu\cdot A_1)\cdot  L^k(\delta_0 \uA, \mbs)  = 0, \ \ \text{for} \ \  k = n-q+1, \ldots, n.
 \end{equation}
This identity is equivalent to saying that \(  \left( \hat \delta_0([\mbu], \uA), \, \mbs\right) = ([\mbu\cdot A_1], A_2, \cdots, A_{r+1},\mbs)\) lies in \(  \theta^{p}_{n,r},\) 
i.e.
\( (\bone\times \partial_0)^{-1}(\theta^{p}_{n,r+1}) =   (\hat\delta_0 \times\bone)^{-1} (\theta^{p}_{n,r}). \) 
The result follows.
\end{proof}

\subsubsection{Proof of the theorem}\hfill

We know that \(  \mathfrak{C}^p_n  \in \sbcxeq{p}{0}{B_\bullet GL_n} \)  and  \( \Theta^q_n \in
\sbcx{q}{0}{\bbp_\bullet(\bbe_n)} \) are cycles  representing classes in their respective higher Chow groups.

\begin{claim}
The cycle \( \Theta^q_n \) represents the \( q\)-th power of the first Chern class of the hyperplane bundle over \( \bbp_\bullet(\bbe_n) \), that is,
\begin{equation}
[ \Theta^q_n ] \ = \ \mbc_1(\scrO_{\bbp_\bullet(\bbe_n)}(1))^q \ \in \ \shcg{q}{0}{\bbp_\bullet(\bbe_n)}. 
\end{equation}
\end{claim}
\begin{proof}[Proof of Claim]

Consider the following diagram, where $pr$ denotes  the projections.
\[ 
\xymatrix{
\bbp^{n-1} & \ar[l]_-{pr_1} \bbp^{n-1}\times GL_n^{\times r}  \times \Delta^r \ar[r]^-{pr_{23}} \ar@/_1pc/[rr]_{\scrL_{n,r}:= L_{n,r}\, \circ\, p_{23}}
& GL_n^{\times r}\times \Delta^r \ar[r]^-{L_{n,r}} & \mat{n}{n} }
\]
The ``$j$-th column morphism''\ $\scrL^j_{n,r} \colon \bbp^{n-1}\times GL_n^{\times r}  \times \Delta^r \to \mat{n}{1} $  can be seen as a section of
$\scrO^{\oplus n}$, for \( j= 1, \ldots, n \).  Therefore, using the quotient map  $\pi\colon \scrO^{\oplus n} \to pr_1^*\scrO_{\bbp^{n-1}}(1)$ we define a section
\[
s_j \in \Gamma( \bbp^{n-1}\times GL_n^{\times r } \times \Delta^r , \pi_1^*(\scrO_{\bbp^{n-1}}(1))
\]
by $s_j = \pi \circ \scrL^j_{n,r}.$ In other words, 
on an \(  R \)-valued point \(  ([\mbu], \uA, \mbt) \)\ the section is given by 
\(
s_j([\mbu], \uA, \mbt) = \la L_{n,r}^j(\uA, \mbt) \ra \in   R^{n}/\mbu^\perp  =   \scrO_{\bbp^{n-1}}(1)_{[\mbu]}(R)
\).

Notice that the  zero set of \(  s_j \) is the subscheme \(  D_{n,r}^1(j) \subset \bbp^{n-1}\times GL_n^{\times r} \times \Delta^r \) whose $R$-valued points are
\begin{equation}
\label{eq:R-val-D}
D_{n,r}^1(j)(R) := \{ ([u], \uA, \mbt) \mid u \cdot L_{n,r}^j(\uA, \mbs) = 0\}.
\end{equation}

Therefore, \(  D_{n,r}^1(j) \)  represents the first Chern class of \(  pr_1^*\scrO_{\bbp^{n-1}}(1) \), for any $j=1,\ldots, n$. By definition, 
\( \theta^{q}_{n,r}  = D_{n,r}^1(n-q+1)\cap \cdots \cap  D_{n,r}^1(n-1)\cap  D_{n,r}^1(n)  \), and a dimension counting shows that this is a proper intersection. A simple  multiplicity calculation now shows that $[\theta^{p}_{n,r}]$
represents the $p$-th power of the first Chern class of $\pi_1^*\scrO_{\bbp^{n-1}}(1)$ over $\bbp^{n-1}\times GL_n^{\times r } \times \Delta^r.$

The  arguments proving Lemma \ref{lem:aux2} show that the elements \(\mathbf{D}^1_n(j) := (D_{n,r}^1(j))_{r\geq 0} \) are special cycles in the complex \( \sbcxeq{1}{0}{\bbp_\bullet(\bbe_n)}\), and it follows from the discussion above that \( [\mathbf{D}^1_n(j)] = \mbc_1(\scrO_{\bbp_\bullet(\bbe_n)}(1)) \in \shcg{1}{0}{\bbp_\bullet(\bbe_n)}.\) Finally,  the fact that  \( \mathbf{D}^1_n(j) \) is a special cycle for j=1, \ldots, n, along with Proposition \ref{prop:special} and  pull-back to the  diagonal, imply that \( [\Theta^q_n] = [\Theta^1_n]\cup \cdots \cup [\Theta^1_n] \). 
\end{proof}

\begin{claim}
Denote \( \xi := \mbc_1(\scrO_{\bbp_\bullet(\bbe_n)}(1)) \in \shcg{1}{0}{\bbp_\bullet(\bbe_n)}\), and let \( \pi \colon \bbp_\bullet(\bbe_n)
\to B_\bullet GL_n \) be the projection. Then the following identity holds:
\[
\sum_{q=0}^n (-1)^q \ \xi^q  \cdot [\frC^{n-q}_n] \ = \ 0.
\]
\end{claim}
\begin{proof}[Proof of Claim]
Let \( \rho \colon \bbp^{n-1}\times \mat{n}{n} \to \mat{n}{n}\) be the projection. 
We have shown in Corollary \ref{cor:sch_int} that \( [\theta^q_n] \centerdot  \rho^*[\bbd^{n-q}] \ = \ [\bbs_{n,n-q}] + [ \bbs_{n,n-q+1}] \), as cycles in \( \bbp^{n-1}\times \mat{n}{n} \). 
On the other hand,  by definition, $C^{n-q}_{n,r} = L_{n,r}^{-1}(\bbd^{n-q}_n)$ and $\theta^{q}_{n,r} = (\bone\times L_{n,r})^{-1}(\theta^q_n)$. Therefore
\begin{align*}
[\theta^{q}_{n,r}]\centerdot \pi^{*}[C^{n-q}_{n,r}] & = (\bone\times L_{n,r})^* \left( [\theta^q_n] \centerdot  \rho^*[\bbd^{n-q}_n]  \right) \\
& = (\bone\times L_{n,r})^*( [\bbs_{n,n-q}] + [\bbs_{n,n-q+1}]).
\end{align*}
It follows that, for each $r\geq 1$, 
\begin{align*}
    \sum_{q=0}^n(-1)^q \ [\theta^{q}_{n,r}] \centerdot \pi^*[{C}^{n-q}_{n,r}] 
     & =  (\bone\times L_{n,r})^*    \left\{ \sum_{q=0}^n\ (-1)^q 
  (\bbs_{n,n-q} + \bbs_{n,n-q+1})
     \right\} = 0,
\end{align*}
since \( [\bbs_{n,0}] = [\bbs_{n,n+1}] = 0 \) by definition.
Now, we invoke Proposition \ref{prop:special} once again to conclude that 
\( \sum_{q=0}^n (-1)^q \ \xi^q \,  \boxtimes  [\frC^{n-q}_n] =  0\). 
Pulling back along the diagonal proves the claim.
\end{proof}

Finally, the compatibility of $\mathfrak{C}^p_n$ with respect to the stabilization maps 
$$ \jmath_{n} \colon B_\bullet GL_n \hookrightarrow B_\bullet GL_{n+1}$$ follows   from 
Proposition \ref{prop:coherent} (or Proposition \ref{prop:cycles2}.\ref{it:2c}). This concludes the proof  of Theorem \ref{thm:main-T}.

\qed

\section{Chern classes in Karoubi-Villamayor \(K\)-theory}
\label{sec:CC-KV}

Let \(X= \spec{R}\in \osch{\Bbbk}\) be an integral affine scheme.
In this section, we give an alternative presentation of  Chern class maps \(\mbc_{p,r} \colon KV_r(R) \to CH^p(X,r)\) from  Karoubi-Villamayor \(K\)-theory to higher Chow groups, utilizing the universal cycles \(\mathfrak{C}^p \)  introduced in the previous section. 
This approach should be contrasted with the usual one via the  homomorphism  
\(
\gamma_{\mathfrak{C}^p} \colon K_r(X) \to CH^p(X;r)
\)
described in \eqref{eq:gamma_c}.  
 
The main difference between the two approaches lies at which stage a moving lemma is applied. The usual approach utilizes the evaluation map
 \(
\ev \colon \mathsf{Spec}(A) \times B_\bullet GL_n(A) \longrightarrow B_\bullet GL_n
\)
 \eqref{eq:eval}  to pull-back the universal Chern classes in \( \shcg{*}{0}{B_\bullet GL(R)}\). This contravariant functoriality requires regularity of \(\spec{A}\) along with a moving lemma for   higher Chow groups. The  approach presented  in this section
uses the homotopy invariance of  \( KV_*(-)\)  and the fact that an element \( \alpha \in KV_r(R) \) can be ``stably represented''   by a morphism with suitable flatness properties, as shown in Lemma \ref{lem:mov} and Proposition \ref{prop:prop}. This is a sort of \emph{moving lemma} for \( KV\)-theory that  holds for  integral \( \Bbbk\)-algebras of finite type. 
 
 \begin{remark} \hfill
  \label{rem:affine}
\begin{enumerate}[i.] 
\item The functors \(KV_j\) are  homotopy invariant, contrary to Quillen's algebraic \(K\)-theory. However, when \( R\) is  regular one has an isomorphism
 \( K_*(R)\cong KV_*(R) \); see \cite[Thm. 11.8]{Weibel_K-book}. In this case, both approaches to Chern classes   coincide.  
\item It is often convenient to write \(KV_r(X)\)  instead of \(KV_r(R)\), where   \(X=\spec{R}\).
\end{enumerate}
 \end{remark}

Let us first recall the definition of \(KV_*(R)\). The simplicial ring
\begin{equation}
    \label{eq:sring}
    R[\Delta^\bullet] \colon\ \  n \mapsto R[\Delta^n] := R[t_0, \ldots, t_n]/\la t_0 + \cdots + t_n - 1\ra,
\end{equation}
gives bisimplicial sets \(  B_\bullet GL_n(R[\Delta^\bullet]) \), \( n\geq 0\); see  Example \ref{exmp:bar}. Define
\( B_\bullet GL (R[\Delta^\bullet]) := \colim{n} B_\bullet GL_n(R[ \Delta^\bullet]) \) and let
\begin{equation}
    \label{eq:KVdiag}
    B_\bullet^\Delta GL(R) := \text{diag} \{ B_\bullet GL(R[\Delta^\bullet])\}
\end{equation}
denote the diagonal simplicial set, with face and degeneracy maps denoted by \( \hat\delta_j\) and \( \hat \sigma_j \), respectively. By definition,
\begin{equation}
    \label{eq:KVKT}
    KV_r(R) := \pi_r \left( B_\bullet^\Delta GL(R), I_\bullet  \right) = \widetilde{B}_r^\Delta GL(R)/\sim_h
\end{equation}
where \(I_\bullet\) is the base-point and
    \[ \widetilde{B}_r^\Delta GL(R) := \{ \mbu \in  B_r^\Delta GL(R)  \mid \hat \delta_j(\mbu) = I_{r-1},\ j=0, \ldots, r \}.\]  
\begin{remark}
\label{rem:htpy}
For completeness, let us recall that the homotopy relation \(\sim_h \)  is generated by \(\mbu \sim_h \mbv \in  \widetilde{B}_r^\Delta GL(R)\) iff there is \(\Gamma \in B_{r+1}^\Delta GL(R)\) such that \(\ \hat \delta_j(\Gamma) = I_r,\)  \( \  j=0, \ldots, r-1,\)\ \(\hat\delta_r(\Gamma) = \mbu\ \) and \(\ \hat \delta_{r+1}(\Gamma) = \mbv; \) see \cite[Defn.~3.6]{May-simplicial}.  Denote this statement by \(\Gamma \colon \mbu \sim_h \mbv.\) 
\end{remark}

Expressed as a morphism, one sees that \(\mbu \) lies in \(\widetilde{B}_r^\Delta GL(R)\) if and only if for every \(j=0, \ldots, r\) the following diagram commutes.
    \begin{equation}
    \label{eq:htpy_comm}
    \xymatrix{
    X  \times \Delta^r \ar[rr]^\mbu  &  & GL^{\times r} \ar[d]^{\delta_j}\\
    X  \times \Delta^{r-1} \ar[dr] \ar[u]^{(1  \times \partial_j)} \ar[rr]_{\hat \delta_j(\mbu)} & &  GL^{\times r-1} \\
    & \spec{\Bbbk} \ar[ur]_{I_{r-1}}.
    }
    \end{equation}

\begin{lemma}
\label{lem:SL}
When \(r\geq 2\) and \( R\) is a domain,   an element \( \mba \in \widetilde{B}_r^\Delta GL_n(R) \) is given by a morphism \( \mba \colon X \times \Delta^r \to SL_n^{\times r} . \)
\end{lemma}
\begin{proof}
Write \( \mba = (a_1, \dots, a_r) \), where \( a_j \colon X\times \Delta^r \to GL_n, \)  \(   j=1, \ldots, r,\) are the coordinate maps. Then \(\det(a_j)\) is an invertible element in the polynomial ring \(  R[\Delta^r] \cong R[y_1, \ldots, y_r]\). Since \( R \) is a domain, it follows that \(\text{det}( a_j) \) lies in \( R^\times \), and the boundary conditions 
on elements in \( \widetilde{B}_r^\Delta GL_n(R) \), with  \( r \geq 2\) , imply that \(\det(a_j) = 1,  \)  for all \( j=1, \ldots, r.\) 
\end{proof}

\subsection{LULU is eventually flat}

Let \( \bbl_{n}, \bbu_{n} \subset  SL_{n} \) be  the closed subgroups of \emph{lower triangular} and \emph{upper triangular unipotent} matrices, respectively, and denote $\bbx_n :=  \bbl_{n}\times \bbu_{n}\times \bbl_{n}\times \bbu_{n}$. Observe that \( \bbx_n \) is an affine space of dimension \( 2n(n-1) \), and  let \( \mu  \colon \bbx_n \to SL_{n}$ be the multiplication map 
\[  \mbA:= (A_1, B_1, A_2, B_2)\mapsto \mu(\mbA) = A_1\cdot B_1\cdot  A_2\cdot  B_2.\]

It follows from \cite{SSV-Gauss} that $\mu$ is surjective, allowing us to use   Nori's equidimensionality lemma (Theorem \ref{thm:equid})  to conclude that whenever \( m\geq n^2-1\) the map
\begin{align}
\label{eq:F}
\mu_m \colon {\bbx_n \times \cdots \times \bbx_n} &\longrightarrow SL_{n} \\
\underline{\mbA}= (\mbA_1, \ldots, \mbA_m) & \longmapsto \mu(\mbA_1) \mu(\mbA_2) \cdots \mu(\mbA_m) \notag
\end{align}
is equidimensional. In this case, $\mu_m$ is actually \emph{faithfully flat} since all schemes involved are smooth over  $\Bbbk$; see \cite[Cor. 14.128]{GortzWed-AG}. 

\begin{remark}
The same result remains true for more general simply-connected Chevalley groups  $G$ over a field, instead of \( SL_n\), for they still admit unitriangular factorizations of length $4$, as shown in  \cite{SSV-Gauss}.
\end{remark}

\begin{definition}
\label{def:contracting}
Given \( n\geq 2\)  let \( \scrF_n \) denote the  product
\[
 \scrF_n  := \underbrace{(\bbx_n \times \cdots \times \bbx_n)}_{n^2-1}  \ \cong \ \bba^{2n(n-1)(n^2-1)},
\]
and define a ``contracting homotopy''  \( \frh_n \colon \scrF_n \times \bba^1 \to SL_n \) by
\begin{equation}
\label{eq:Hrn}
{\frh}_n (\underline{\mbA} ; t)  \ := \  \mu_{n^2-1}(\underline I + t( \underline\mbA- \underline I) ),
\end{equation} 
where  \( \underline I \in \bbx_n^{n^2-1} \) is the identity element for the product group-scheme structure, and  addition is defined as the componentwise addition of matrices. 
\end{definition}

\begin{proposition}
\label{prop:H}
 The following properties hold.

\begin{enumerate}[i.]
\item \label{it:H1} \( \frh_n(\underline \mbA;0) = {I}\) for all \( \underline \mbA \in \scrF_n\), where \( {I}\) is the identity element in \( SL_n\).
\item \label{it:H2}  Given \( y \in \bba^1- \{ 0 \} \), the restriction \(  \frh_{n|y} \colon \scrF_{n|\Bbbk(y)} \to SL_{n|\Bbbk(y)}\) is faithfully flat.
\item \label{it:H3}  The natural inclusion \( \iota_n \colon \scrF_n \hookrightarrow \scrF_{n+1} \)  induced by \( \jmath_n \colon SL_n \hookrightarrow SL_{n+1} \), embeds \( \scrF_n \) as an affine subspace of \( \scrF_{n+1}\), yielding a commutative diagram.
\[
\begin{tikzcd}
\scrF_n \times \bba^1   \ar[rr, "\frh_n"] \ar[d, hook, "\iota_n \times \bone"'] &  & SL_n \ar[d, hook, "\jmath_n"] \\
\scrF_{n+1} \times \bba^1 \  \ar[rr, "\frh_{n+1}"] &  & SL_{n+1} 
\end{tikzcd}
\]
\end{enumerate}
\end{proposition}

\begin{proof}
The first and third assertions follow directly from the definitions. The second assertion is a corollary of Theorem \ref{thm:equid}, \eqref{eq:F} and \eqref{eq:Hrn}, since \( \underline\mbA  \mapsto \underline I + t\left( \underline\mbA - \underline I\right) \) defines an  automorphism of \( \bbx^{n^2-1}_{n|\Bbbk(t)} \) over \( \Bbbk(t)\).
\end{proof}

\subsection{A stable moving lemma for \( KV \)-theory}

We now    describe a canonical way to represent a class \( \alpha \in KV_r(X)\), \( r\geq 2\), 
by an element \(\lambda_\mba \colon (X\times \bba^N)\times \Delta^r \to SL_n^{\times r}\) whose restriction to \( (X\times \bba^N)\times (\Delta^r - \partial \Delta^r) \) is faithfully flat. Here \( N=\phi(n,r) \) is an explicit function of \( n\) and \( r\),  when  \(\alpha\) has a representative \( \mba \colon X\times \Delta^r \to SL_n^{\times r}\).
The various homotopies in the construction   carry similar flatness properties that yield corresponding homology relations amongst cycles in the higher Chow complexes. 

\begin{definition}
\label{def:lambda}
For any \( r\geq 0\) define \( \rho \colon \Delta^r \to \bba^1\)   by  \( \rho(\mbt) = t_0 t_1 \cdots t_r\) and, for \( \ell =0, \ldots, r \),  define \( \rho_\ell(\mbt):= t_0 \cdots \widehat{t}_\ell \cdots t_r\). Denote \( \scrF_{n,r} := \scrF_n^{\times (r+1)} \).
Given  \( \mba \colon X \times \Delta^r \to SL_n^{\times r} \), with coordinate functions \( \mba = (\mba_1, \ldots, \mba_r),\) define
\begin{align}
\label{eq:R}
\lambda_\mba \colon X \times \scrF_{n,r} \times \Delta^r & \longrightarrow SL_n^{\times r} \\
(x; \mbA_0|\mbA_1| \cdots  |\mbA_{r}; \mbt) & \longmapsto  
\left(
\mba_1(x,\mbt)\cdot \frh_n(\mbA_1 ,  \rho(\mbt)) | \cdots | \mba_r(x,\mbt)\cdot  \frh_n(\mbA_r ,  \rho(\mbt) )
\right) \notag
\end{align}
Given \( 0\leq i_1 < \cdots < i_k \leq r \),  let \( Z_{i_1, \ldots, i_k} \subset \partial \Delta^r\) be the subscheme defined by the monomial ideal
\( \la \rho_{i_1}(\mbt), \ldots, \rho_{i_k}(\mbt) \ra \).
\end{definition}
\begin{remark}
Note that the definition of \( \lambda_\mba\) does not involve the first component \( \mbA_{0}\) of the product \( \scrF_{n,r}\). This extra factor will be used to create  ``enough room''  for suitably flat homotopies in subsequent arguments. 
\end{remark}

\begin{lemma}
\label{lem:mov}
Let  \( \mba \colon X \times \Delta^r \to SL_n^{\times r} \)  \ be as above. Then
for all \( \mbt \in \Delta^n - \partial \Delta^n\)  the morphism
\( \lambda_{\mba |\mbt} \colon (X\times \scrF_{n,r})_{\Bbbk(\mbt)} \to SL_{n|\Bbbk(\mbt)}^{\times r} \)
is faithfully flat.
\end{lemma}
\begin{proof}
The result follows directly  from Lemma \ref{lem:flat}, Proposition \ref{prop:H} and \cite[Cor. 14.128]{GortzWed-AG}.
\end{proof}

Let  \( R[\mbx_{n,r}]\) be  a polynomial ring over \( R\) such that \( X\times \scrF_{n,r}\cong \spec{R[\mbx_{n,r}]}=\bba^{\phi(n,r)}_R \),
where \( \phi(n,r) = 2(r+1)n(n-1)(n^2-1)\).

\begin{proposition}
\label{prop:prop}
Let \( X \) be an integral affine scheme in \(\osch{\Bbbk}\). Fix \( r \geq 2\) and let \( \pi \colon X \times \scrF_{n,r}   \to X    \) denote the projection. 
Given \( \mbc\colon X\times \Delta^r \to SL_n^{\times r} \), one has
\begin{equation}
\label{eq:restr}
\lambda_{\mbc}|_{X\times \scrF_{n,r}\times \partial\Delta^r} \ = \ 
\mbc\circ( \pi \times \bone)|_{X\times \scrF_{n,r}\times \partial\Delta^r};
\end{equation}
see \eqref{eq:R}. Furthermore, if \( \mba, \mbb \in \widetilde B_r^\Delta GL_n(R) \)  the following holds.
\begin{enumerate}[\rm i.]
\item  \label{it:1} \( \lambda_\mba, \lambda_\mbb \in \widetilde B_r^\Delta SL_n(R[\mbx_{n,r}]). \) (See Lemma \ref{lem:SL}. )
\item \label{it:2}  \( \lambda_\mba \  \sim_h \, \pi^*\mba :=  \mba \circ (\pi\times \bone)  \). (See Remark \ref{rem:htpy}.)
\item \label{it:3} Given   \( \Gamma \colon \mba \sim_h \mbb\ \) there is 
\( \ \tilde\Gamma \colon \lambda_\mba \sim_h \lambda_\mbb \)
such that the restrictions
\[
\tilde \Gamma_{\mbt} \colon \left(X \times \scrF_{n,r}\right)|_{\Bbbk(\mbt)} \to SL_{n|\Bbbk(\mbt)}^{\times r}
\]
are faithfully flat for all \(  \mbt \in \Delta^n -  Z_{r,r+1} \),
where \( Z_{r,r+1} \subset \partial \Delta^{r+1} \) is the closed subscheme introduced in Definition \ref{def:lambda}.
\item  \label{it:4}  If  \( \Gamma \colon X\times \Delta^{r+1} \to SL_n^{r+1} \) 
 `fills of the horn' \( (I, \ldots, I, \mba, -, \mbb) \), that is,
\[ \Gamma|_{X\times \partial_j \Delta^{r+1}} = 
\begin{cases}
I &, \text{ if } j=0,\ldots, r-2 \\
\mba &, \text{ if } j = r-1 \\
\mbb &, \text{ if } j = r+1.
\end{cases}.\]
then there is  \( \tilde\Gamma \colon (X\times \scrF_{n,r})\times \Delta^{r+1} \to SL_n^{r+1} \) satisfying:
\begin{enumerate}[a.]
\item \( \tilde \Gamma\) fills the horn \( (I, \ldots, I, \lambda_\mba, - , \lambda_\mbb) \).
\item If \( \mbc := \Gamma|_{X\times \partial_r\Delta^{r+1}}\)  \ then \ 
\( \tilde \Gamma_{|X\times \scrF_{n,r} \times \partial_r\Delta^{r+1}} \ = \  \lambda_\mbc\).
\item For all 
\( \mbt \in \Delta^n -  Z_{r-1, r, r+1} \) the restriction 
\( \tilde \Gamma_{\mbt} \colon (X \times \scrF^{(r+1)}_n)_{\Bbbk(\mbt)} \to SL_{n|\Bbbk(\mbt)}^{\times r} \) is faithfully flat.
\end{enumerate}
\end{enumerate}
\end{proposition}
\begin{proof}
See Appendix \ref{subsec:pf_prop:prop}.
\end{proof}

The constructions above give the following ``stable moving lemma'' for \( KV\)-theory.

\begin{corollary}
\label{cor:KV-moving}
Let  \( \alpha \in KV_r(X), r \geq 2, \)  be represented by  \( \mba \colon X \times \Delta^r \to SL_n^{\times r} \). Then, there exists a canonical 
\( \lambda_\mba \colon X\times \scrF_{n,r} \times \Delta^r \to SL_n^{\times r} \) satisfying the following.
\begin{enumerate}[i.]
\item The restriction of \( \lambda_\mba\) to \( X\times \scrF_{n,r} \times (\Delta^r - \partial\Delta^r) \) is faithfully flat.
\item \( \lambda_\mba\) represents a class \( \lambda_\alpha \in KV_r(X\times \scrF_{n,r})\) corresponding to \( \alpha \) under the isomorphism
\( \rho^*\colon KV_r(X) \xrightarrow{\ \cong \ } KV_r(X\times \scrF_{n,r}) \) induced by the projection \( \rho \colon X\times \scrF_{n,r} \to X \).
\end{enumerate}
\end{corollary}

\subsection{The Chern class homomorphims $ \mbc_{p,r} \colon KV_r(R) \to \shcg{p}{r}{X}$}

Consider  an integral \( \Bbbk\)-algebra \( R \) of finite type, and denote \( X = \spec{R} \).

\subsubsection{\bf The case  $\mathbf{r\geq 2}$}
\label{subsubsec:rgeq2}
Consider the commutative diagram 
\[
\begin{tikzcd}
GL_n^{\times r} \times \Delta^r  \ar[rr, "(\pi_k \circ L)\times 1"]   & & \mat{n}{k} \times \Delta^r \\
SL_n^{\times r} \times \Delta^r  \ar[hook, "\jmath"]{u}\ar[rru, "J_{n,k}"'] & & 
\end{tikzcd}
\]
It follows from  Lemma \ref{lem:L-smooth}, that both \( (\pi_k\circ L) \times 1\) and \( J_{n,k} \) are smooth when restricted to \( GL_n^{\times r} \times ({\Delta}^r-\partial \Delta^r)  \) and 
\( SL_n^{\times r} \times ({\Delta}^r-\partial \Delta^r )\), respectively. Therefore, 
the subschemes \( C^p_{n,r}\) and \( SL_n^{\times r} \times \Delta^r\)   intersect properly in  \(GL_n^{\times r} \times \Delta^r\),  and their intersection is still a Cohen-Macaulay, irreducible subscheme  of \(  SL_n^{\times r} \times \Delta^r\)  that we denote \( \jmath^*C^p_{n,r} \).

Given   \( f \colon Y \times \Delta^r \to GL_n^{\times r} \),  denote by   \( \hat f \colon Y \times \Delta^r \to GL_n^{\times r} \times \Delta^r\)   the morphism   
\begin{equation}
\label{eq:fhat}
\hat f \colon  (x, \mbt) \mapsto ( f(x,\mbt),\, \mbt).
\end{equation}
\begin{proposition}
\label{lem:flat1}
Let \( R \) be an integral \( \Bbbk\)-algebra, and fix \( 1\leq p \leq n\) and \( r\geq 2\). 
The assignment that sends 
\( \mba \in \widetilde{B}_r^\Delta GL_n(R) \)
to the flat pull-back
\( \hat \lambda_\mba^*\left[ \jmath^*C^p_{n,r}\right] \)
gives a well-defined function
\begin{equation*}
\tilde \gamma_n^p \colon \widetilde{B}_r^\Delta GL_n(R)  \longrightarrow  \widetilde{\scrZ}^p_\text{eq}(X\times \scrF_{n,r};r) ,
\end{equation*} 
where \(\widetilde{\scrZ}^p_\text{eq}(X\times \scrF_{n,r};r)  := \ker\{ \sbcxeq{p}{r }{X\times \scrF_{n,r}} \xrightarrow{\partial}\sbcxeq{p}{r-1}{X\times\scrF_{n,r}} \}.\)
This assignment induces a homomorphism
 \begin{equation}
 \label{eq:fun_HCG}
  \gamma_n^p \colon \pi_r \left( B_\bullet^\Delta GL_n(R) \right)   \to \shcg{p}{r}{X\times \scrF_{n,r}},
  \end{equation}
 compatible with stabilizations. In other words, the following diagram commutes:
  \[
  \begin{tikzcd}
  \pi_r \left( B_\bullet^\Delta GL_n(R) \right) \ar[dd, "\jmath_n"']   \ar[r, "\gamma_n^p"] &  \shcg{p}{r}{X\times \scrF^r_n} \ar[dd, "pr^*"', "\cong"] \ar[dr, "(pr^*)^{-1}", "\cong"'] & \\
& & \shcg{p}{r}{X}. \\
    \pi_r \left( B_\bullet^\Delta GL_{n+1}(R) \right)   \ar[r, "\gamma_{n+1}^p"'] &  \shcg{p}{r}{X\times \scrF^r_{n+1}}\ar[ur, "(pr^*)^{-1}"', "\cong"]  & \\
  \end{tikzcd}
  \]
\end{proposition}
\begin{proof}
Since \(\jmath^* C^p_{n,r}\) is equidimensional and dominant over \( \Delta^r\) (see  Lemma~\ref{lem:L-smooth}.\ref{it:L-b}), it follows from 
Lemma \ref{lem:mov} that  \( \hat\lambda_\mba^*[\jmath^* C^p_{n,r}] \) lies in \( \sbcxeq{p}{r}{X\times\scrF_{n,r}}. \) 
On the other hand,  Proposition \ref{prop:prop}.\ref{it:1} and   the definition of  \( C^p_{n,r}\) imply that whenever  \( \mba\) lies in \( \widetilde{B}_r^\Delta GL_n(R)  \)  the cycle
\( \hat\lambda_\mba^*[\jmath^* C^p_{n,r}] \) does not intersect \( X \times \scrF_{n,r} \times \partial \Delta^r \), and hence it lies in 
\( \widetilde{\scrZ}^p_\text{eq}(X\times \scrF_{n,r};r) \).

Finally,  Proposition \ref{prop:prop} shows that \( \gamma_n^p \) is a homomorphism, whose compatibility with stabilizations follows from Proposition \ref{prop:H}.\ref{it:H3}, Theorem \ref{thm:main-T} and   definitions. 
\end{proof}

In summary,  the homomorphisms \( (pr^*)^{-1} \circ \gamma^p_n \colon  \pi_r \left( B_\bullet^\Delta GL_n(R) \right) \to \shcg{p}{r}{X}\)  induce a well-defined \emph{Chern class} homomorphism
\begin{equation}
\mbc_{p,r} \colon KV_r(R) = \colim{n}\pi_n\left(B_\bullet^\Delta GL_n(R)\right)\ \  \longrightarrow\ \  \shcg{p}{r}{X}.
\end{equation}

\begin{corollary}
\label{cor:chernKV}
For \( r \geq 2\) and \(  p\geq 1\), the  homomorphism
\(
\mbc_{p,r}  \)
is a natural transformation between the functors \(KV_r(-)\) and \(\shcg{p}{r}{\spec{-}} \) in the category of regular \(\Bbbk\)-algebras.
\end{corollary}

\subsubsection{{\bf The case  $\mathbf{r=1} $} (part I): $R^\times \to \shcg{1}{1}{X}$.}\hfill
\label{subsec:baby}
\smallskip 

Recall that 
\(
K_1(R) := GL(R)/E(R) \cong R^\times \oplus SK_1(R),
\)
where \( E(R) = \colim{n}E_n(R) \) is the subgroup of elementary matrices and \( SK_1(R) = SL(R)/E(R)\).
From now on, assume that \( R\) is a \( \Bbbk\)-algebra of finite type and, for \(n\geq~1\), 
denote   \(E_n := E_n(R[x][\Delta^1]))\subset GL_n(R[x][\Delta^1])\).
Given \(\alpha~\in~R^\times\) define 
\begin{equation}
\label{eq:Aalpha}
A_\alpha := \begin{pmatrix} x t_0 & x^2 t_0t_1 - \alpha \\ 1 & xt_1 \end{pmatrix} \in GL_2(R[x][\Delta^1]). 
\end{equation}

\begin{lemma}
\label{lem:comm}
Given \(\alpha, \beta \in R^\times\), define 
\[
U_\beta:= \begin{pmatrix}
   1 & \beta(x-1)- xt_1 \\
  0 & 1
\end{pmatrix}, \quad
V_\beta:= \begin{pmatrix}
 1 & 0 \\
   \beta^{-1} & 1
\end{pmatrix}
\quad \text{and} \quad 
W_\beta :=\begin{pmatrix}
  1 & xt_1 - \beta \\
 0 & 1
\end{pmatrix}.
\]
Then
\(
A_{\alpha\beta}^{-1} A_\alpha A_\beta \ = \ U_\beta V_\beta W_\beta.
\)
In particular, 
\begin{enumerate}[a.]
\item \(A_{\alpha \beta}^{-1}A_\alpha A_\beta \in E_2\).
\item \(A_1 \in E_2\)
\item \(A_{\alpha^{-1}} A_\alpha \in E_2 \)
\end{enumerate}
\end{lemma}
\begin{proof}
The identity \( A_{\alpha\beta}^{-1} A_\alpha A_\beta \ = \ U_\beta V_\beta W_\beta \)  follows from a straightforward calculation, using the fact that \( t_0+t_1 = 1 \), and assertion (a) is a direct consequence of this identity. Assertion (b) follows  by setting \(\alpha =\beta=1\) in the identity, while
assertion (c) follows by setting \(\beta = \alpha^{-1}\).
\end{proof}

Denote \(E:= E(R[x][\Delta^1]) = \colim{n} E_n(R[x][\Delta^1])\) and let \(\eta \colon R^\times \to K_1(R[x][\Delta^1])  \) be the composition
\begin{align}
    \label{eq:eta}
\eta \colon R^\times & \longrightarrow  GL_2(R[x][\Delta^1])/ E_2(R[x][\Delta^1]) \to GL(R[\mbx][\Delta^1])/ E(R[x][\Delta^1])   \\
\alpha & \longmapsto A_\alpha. \notag
\end{align}
It follows from the lemma above that \(\eta\) is a homomorphism and the composition
\(
R^\times \xrightarrow{ \ \eta \ }  K_1(R[ x][\Delta^1]) \xrightarrow{\ \det \ } R^\times
\)
is the identity. In particular, \(\eta\) is an isomorphism when \(R\) is a field, or a regular local ring.


In order to explicitly describe the composition  
\begin{equation}
\label{eq:compo}
R^\times \xrightarrow{\eta} K_1(R[x][\Delta^1])\to  KV_1(R[x])    
\xrightarrow{\mbc_{1,1}} \shcg{1}{1}{X\times \bba^1}\cong \shcg{1}{1}{X},
\end{equation}
consider  \(A_\alpha\) as a morphism
\(A_\alpha  \colon (X\times \bba^1) \times \Delta^1 \longrightarrow GL_2 \). By definition, \(\mbc_{1,1}(\eta(\alpha)) = [\hat A_\alpha^*(C^1_{2,1} )] \in \shcg{1}{1}{X\times \bba^1}\). This requires the verification that the divisor \( \hat A_\alpha^*(C^1_{2,1} )\)  is equidimensional over \( \Delta^1\).
Indeed,  from the definitions of \(C^1_{2,1}\) and \(A_\alpha\) it follows that this divisor is given by
\begin{align}
\label{eq:palpha}
p_\alpha( t_0, t_1;x):= \det( t_0 I + t_1 A_\alpha)  & = \det \begin{pmatrix}
t_0 + xt_0t_1 & t_1(x^2t_0t_1 - \alpha) \notag \\
t_1 & t_0+ xt_1^2
\end{pmatrix}  \notag \\
& = t_0^2 + x t_0t_1 + \alpha t_1^2\ \in \ R[\Delta^1][x],
\end{align}
since \( t_0 + t_1 = 1 \).  In other words, 
\begin{equation}
\label{eq:Rstar}
\mbc_{1,1}\circ \eta(\alpha)= [ \text{div}( t_0^2 + x t_0 t_1 + \alpha t_1^2)] \in \shcg{1}{1}{X \times \bba^1}.
\end{equation}

\begin{remark}
\label{rem:unit}
Observe that  \( \la p_\alpha(t_0, t_1; x) \ra + \la t_0 t_1 \ra = \la 1 \ra\) whenever \( \alpha\) is a unit in \( R \).
\end{remark}

\subsubsection{{\bf The case $\mathbf{r=1}$} (part II):  \({ SK_1(R) \to \shcg{1}{1}{X}}\) }\hfill

The arguments   here are similar to the ones used in \S \ref{subsubsec:rgeq2}.\ 
Given \( \mba \colon X  \to SL_n\),  just as in  Definition \ref{def:lambda}   define
\begin{align*}
\lambda_\mba \colon X \times \scrF_{n,1} \times \Delta^1 & \longrightarrow SL_n \\
(x;\mbA_0|\mbA_1; \mbs) & \longmapsto \mba(x)\, \frh_n(\mbA_1; s_0s_1),
\end{align*}
and  denote  \(D_\mba := \hat \lambda_\mba^*[ C^1_{n,1}] \in \sbcxeq{1}{1}{X\times \scrF_{n,1}} \); see \eqref{eq:fhat} and \eqref{eq:Cs}.

\begin{lemma}
Let $\mba, \mbb \in SL_n(R)$ be given, and denote \( X = \spec{R}\).
\begin{enumerate}[i.]
\item The relation \( D_{\mba\mbb}  \sim  D_\mba + D_\mbb \) holds in the higher Chow complex. Therefore, one obtains a well-defined homomorphism 
\begin{align*}
 \rho_n \colon SL_n(R) & \to \shcg{1}{1}{X\times \scrF_{n,1}} \cong \shcg{1}{1}{X} \\
 \mba & \mapsto \left[ D_\mba \right].
 \end{align*}
\item The maps \( \rho_n\) are compatible with stabilization and thus give a homomorphism \( \rho\colon SL(R) \to \shcg{1}{1}{X}\).
\item The group of elementary matrices \( E(R)\subset SL(R)\) lies in the kernel of \( \rho \).
\end{enumerate}
\end{lemma}
\begin{proof}
Define 
\( 
g_{\mba,\mbb}  \colon X \times \scrF_{n,1} \times \Delta^2  \to SL^{\times 2}_n \)
by
\begin{equation*}
g_{\mba,\mbb} (x; \mbA_0  | \mbA_1; \mbt)  := \left( \mba(x) \frh_n(\mbA_0; t_0t_1t_2) \frh_n(\mbA_1; t_0t_1) \  | \   \mbb(x) \frh_n(\mbA_1; t_0t_2+t_1t_2)\right).
\end{equation*}
This is an element in \( B^\Delta_2SL_n(R[\mbx_{n,1}]) \) satisfying
\begin{align*}
(d_0 g_{\mba,\mbb}) (x;\mbA_0  | \mbA_1; \mbs) & = 
\mbb(x)\frh_n(\mbA_1;s_0s_1) = \lambda_\mbb(x,\mbA_0|\mbA_1;\mbs) \\
(d_1 g_{\mba,\mbb}) (x; \mbA_0  | \mbA_1; \mbt) & =
 \mba(x) \frh_n(\mbA_0; 0)\frh_n(\mbA_1; 0)   \mbb(x) \frh_n(\mbA_1;s_0s_1) \\
 & = \mba(x) \mbb(x) \frh_n(\mbA_1;s_0s_1)  =
   \lambda_{\mba\mbb}(x,\mbA_0|\mbA_1;\mbs)\\
(d_2 g_{\mba,\mbb}) (x;\mbA_0  | \mbA_1; \mbt)& = \mba(x) \frh_n(\mbA_0; 0)\frh_n(\mbA_1; s_0s_1)  =
    \lambda_{\mba}(x,\mbA_0|\mbA_1;\mbs).
\end{align*}

As in the proof of Proposition \ref{prop:prop}, it follows if if \(Z := V\la t_0t_1t_2,  t_0t_1, t_0t_2+t_1t_2  \ra \subset \partial \Delta^2\) then for all \( \mbt \in \Delta^2 -Z \) the restriction
\( g_{\mba,\mbb| \mbt} \) is a faithfully flat map. 
Hence, one can define
\( H_{\mba,\mbb}:= \hat g_{\mba,\mbb}^*[C^1_{n,2}]    \in \sbcxeq{1}{2}{X\times \scrF_{n,1}}\).
It follows from the identities above and definitions that \(\partial H_{\mba,\mbb} = D_\mbb - D_{\mba\mbb}+D_\mba\).
Therefore, \( \rho \colon SL_n(R) \to \shcg{1}{1}{X} \) is a group homomorphism.

The second assertion  follows from Proposition \ref{prop:H} and  Proposition \ref{prop:cycles2}.\ref{it:2c}.

For the last assertion, consider \( \alpha \in E_n(R) \subset SL_n(R) \), and let  \(  \mathsf{p}_\alpha \colon X \times \Delta^1 \to  SL_n  \ \in \ B^\Delta_1SL_n(R) \) be any \emph{path}
satisfying
\( \mathsf{p}_\alpha(x;0,1)  = I \) and \( \mathsf{p}_\alpha(x;1,0) = \alpha(x)\). Then, define \( F_\alpha \colon X \times \scrF_{n,1} \times \Delta^2 \to SL_n^{\times 2} \) by
\begin{multline*}
F_\alpha(x,\mbA,\mbt) = \\
\left(
\mathsf{p}_\alpha (x;t_0+t_1,t_2) \frh_{n}(\mbA_0;t_0 t_1 t_2) \frh_{n}(\mbA_1;t_0 t_1)  \mid
\mathsf{p}^{-1}_{\alpha}(x;t_0,t_1 + t_2)
\frh_n(\mbA_1; t_0 t_1 ) 
\right).
\end{multline*}
Once again, if \(Z := V\la t_0t_1t_2,  t_0t_1 \ra \subset \partial \Delta^2\) then for all \( \mbt \in \Delta^2 -Z \) the restriction
\(F_{\alpha | \mbt} \) is  faithfully flat map,  as in Proposition \ref{prop:prop}.  A quick inspection of the definitions shows that 
\begin{align*}
(d_0 F_{\alpha})(x;\mbA_0 | \mbA_1;\mbs) & =  I \\
(d_1 F_{\alpha} )(x;\mbA_0 | \mbA_1;\mbs)&  = I\\
(d_2 F_{\alpha})(x;\mbA_0 | \mbA_1;\mbs) & = 
    \lambda_{\alpha}(x;\mbA_0|\mbA_1;\mbs).
\end{align*}
Now, define
\( H_{\alpha}:= \hat F_{\alpha}^*[C^1_{n,2}]    \in \sbcxeq{1}{2}{X\times \scrF_{n,1}} \) 
and conclude that \( \partial H_\alpha = D_\alpha\).
It follows that  we have a well-defined  homomorphism
\begin{equation}
\label{eq:rho}
\rho \colon SK_1(R) := SL(R)/E(R) \longrightarrow \shcg{1}{1}{X}.
\end{equation}
\end{proof}

Combining \eqref{eq:rho} and  \eqref{eq:Rstar} gives the desired homomorphism
\begin{equation}
\mbc_{1,1} \colon K_1(R) \cong R^\times \oplus SK_1(R) \longrightarrow \shcg{1}{1}{X}.
\end{equation}

\begin{remark}
When \( X\) is a complex smooth quasi-projective variety \( X\),     
the equidimensional complex  can be used to construct explicit    regulator maps \[ \mathsf{Reg}\colon \shcg{n}{p}{X} \to \intdc{n}{p}{X}; \] see \cite{REMAHI-CHOGRO}. For \(  X = \spec{R} \) affine, if  \( D = \text{div}( f)  \subset X\times \bba^1 \times \Delta^1 \) is an element in 
\( \sbcxeq{1}{1}{X\times \bba^1}\) representing a class in \( \shcg{1}{1}{X\times \bba^1} \cong \shcg{1}{1}{X}\cong \scrO^\times(X)=R^\times\) which is given by a rational function \( f (x;s_0,s_1) \), it is  shown in \cite[\S 5.1]{REMAHI-CHOGRO}  that \( \mathsf{Reg}(D) = \frac{f(x;0,1)}{f(x;1,0)} \in R^\times.\)  The equidimensionality property makes this a direct and simple observation.

Now, consider the composition  \( R^\times \xrightarrow{\eta} K_1(R) \xrightarrow{\mbc_{1,1}} \shcg{1}{1}{X\times \bba^1} \)  in   \eqref{eq:compo}, where for \( \alpha\in R^\times\) the first Chern class \( \mbc_{1,1}(\eta(\alpha)) \in \shcg{1}{1}{X\times \bba^1}\) is represented by \( \text{div}(p_\alpha) \), where \( p_\alpha = s_0^2 + s_0s_1 x + s_1^2 \alpha\). From the previous paragraph one obtains
\[
\mathsf{Reg}(\mbc_{1,1}(\eta(\alpha))) =\frac{p_\alpha(x;0,1)}{p_\alpha(x;1,0)} = \alpha \in R^\times,
\]
as expected.
\end{remark}

\subsection{A Steinberg symbol}

In this section,   assume that \( X=\spec{R} \) is smooth over \( \Bbbk\).
Given \( \alpha\in R\), consider \( p_\alpha (s_0, s_1;y) = s_0^2 + s_0 s_1 y + s_1^2 \alpha \)  as in \eqref{eq:palpha}, and let 
 \( D_\alpha = \text{div}(p_\alpha)\) be the corresponding divisor in \( X\times \bba^1 \times \Delta^1 \).
We have seen  that if  \( \alpha\) lies in \( R^\times \) then \( D_\alpha\) represents the class in \( \shcg{1}{1}{X\times \bba^1}\) corresponding to \( \alpha \) under the isomorphism \( \shcg{1}{1}{X\times \bba^1}\cong R^\times \).
Consider  \( \scrA := R[\mbs, \mbt]/\la s_0+s_1-1,\, t_0 + t_1 -1 \ra = R[\Delta^1\times \Delta^1] \)
and denote
 \begin{equation}
\label{eq:pq}
\mathfrak{p}_\alpha  =  s_0^2 + s_0s_1 y_1 + \alpha y_1^2  \quad \text{and}  \quad 
\mathfrak{q}_\beta   = t_0^2 + t_0t_1 y_2 +  \beta t_1^2 \ \in  \ \scrA [y_1,y_2].
\end{equation}
Now, let \(\pi_1, \pi_2 \colon X \times \bba^2 \times \Delta^1 \times \Delta^1 \to X \times \bba^1 \times \Delta^1 \) denote the projections 
\[
\begin{tikzcd}
   & \ar[dl, mapsto, "\pi_1"'] (x; y_1, y_2; s_0, s_1; t_0, t_1)  \ar[rd, mapsto, "\pi_2"] &   \\
(x; y_1; s_0, s_1)  &   &  (x; y_2; t_0, t_1).
\end{tikzcd}
\]
Given \( \alpha, \beta \in R^\times, \ \) the intersection of divisors
 \(  \pi_1^* D_\alpha \,  \centerdot  \, \pi_2^*D_\beta 
\ = \ \text{div}(\frp_\alpha) \centerdot  \text{div}(\frq_\beta)\)
gives a cycle of codimension \( 2 \) \ in\  \( X \times \bba^2 \times \Delta^1 \times \Delta^1 \) and, using  the Eilenberg-Zilber map
 \( \psi_{1,1} \colon \Delta^2 \to \Delta^1 \times \Delta^{1,1} \), 
 one can  define 
%
\begin{align}
\Gamma_{\alpha, \beta} & :=  \psi_{1,1}^*\left( \text{div}(\frp_\alpha) \centerdot  \text{div}(\frq_\beta) \right).
 \end{align}
Although the properties listed in the next result follow directly from the results proven thus far, we outline a direct proof, as it might be of independent interest.

\begin{proposition}
Given \( \alpha, \beta, \gamma\, \in R^\times \), the following holds.
\begin{enumerate}[a.]
\item  \( D_{\alpha\beta} \sim D_{\alpha}+ D_{\beta} \) and \( D_1 \sim 0\).
\item \(\Gamma_{\alpha,\beta} \) is a cycle in \( \sbcxeq{2}{2}{X\times \bba^2}\).
\item \( \Gamma_{\alpha\beta, \gamma} \sim \Gamma_{\alpha, \gamma} + \Gamma_{\beta, \gamma} \).
\item \(  \Gamma_{\alpha, -\alpha} \sim 0\);
\item If \( \alpha\) and \( 1-\alpha\) are in \( R^\times \) then \( \Gamma_{\alpha, - \alpha }  \sim 0 \);
\end{enumerate}
\end{proposition}

\begin{proof}
Define  \( f(s_0, s_1, y, T) = 1 + s_0s_1(y-2)T \in R[\Delta^1][y,T] \) and observe that 
\( F := \text{div} (f) \) is a divisor in \( (X \times \bba^1) \times \Delta^1 \times \bba^1 \) which is equidimensional and dominant over \( \Delta^1 \times \bba^1\). Since 
\( F_{|s_0=0} = F_{|s_1=0}=F_{|T=0}=0 \) while \( F_{|T=1} = D_1 \),  it follows from Remark \ref{rem:prism} that
\begin{equation}
\label{eq:D1}
D_1 \ \sim \ 0.
\end{equation}

Now, consider
\[
h_{\alpha, \beta} := 
s_0^2 + 2 s_0^3s_1 y + (y^2 + C_{\alpha, \beta}(T) )s_0^2 s_1^2 + y C_{\alpha, \beta}(T) s_0 s_1^3 + \alpha \beta s_1^4\ 
 \in\  R[\Delta^1][y,T],
\]
where \(C_{\alpha, \beta}(T) = (\alpha\beta - \alpha-\beta+1)T + \alpha + \beta\), 
and define \( H_{\alpha, \beta} := \text{div}(h_{\alpha, \beta} ) \) in \(    (X\times \bba^1)\times \Delta^1 \times \bba^1 \).
Note that:
\begin{enumerate}[a.]
\item \( h_{\alpha,\beta}(s_0,s_1,y,0) = p_{\alpha}(x,s_0,s_1)p_{\beta}(x,s_0,s_1) \);
\item \(  h_{\alpha,\beta}(s_0,s_1,y,1) = p_{\alpha\beta}(x,s_0,s_1) p_{1}(x,s_0,s_1)\);
\item \(  h_{\alpha,\beta}(s_0,s_1,y,T) \) is a polyomial of degree \( 2\) in the variable \( y \), whose leading coefficient is  \( s^2_0s^2_1 \);
%
\item \( \la h_{\alpha, \beta} \ra \ + \ \la s_0s_1 \ra = \la 1 \ra \).
\end{enumerate}
These properties suffice to show that \( H_{\alpha, \beta} \) is equidimensional and dominant over \( \Delta^1\times \bba^1 \), and that 
\( H_{\alpha, \beta|_{T=0}} = \text{div}(p_{\alpha}p_\beta)    = D_\alpha + D_\beta \) and  \( H_{\alpha, \beta|_{T=1}} = \text{div}(p_{\alpha \beta}p_{1}) = D_{\alpha\beta}+D_1. \) Since, \( H_{\alpha, \beta|_{s_0=0}} = 0 \) and  \( H_{\alpha, \beta|_{s_1=0}} = 0\), the first assertion then follows from Remark \ref{rem:prism} and \eqref{eq:D1}.

To prove the second statement, 
let  \( \scrJ_{\alpha, \beta}\subset \scrA[y_1,y_2] \) be ideal \( \scrJ_{\alpha, \beta} :=\la \frp_\alpha, \frq_\beta \ra\). A direct calculation shows that 
\( \la \frp_\alpha \ra + \la s_0s_1\ra = \la 1 \ra = \la \frq_\beta \ra + \la t_0t_1 \ra \), and hence \( s_0, s_1, t_0, t_1 \) are all units in \( \scrA[y_1, y_2 ]/\scrJ_{\alpha, \beta} \). 
Therefore,
\begin{equation}
\label{eq:inter1}
V(\scrJ_{\alpha, \beta}) \cap  (X\times \bba^2 \times \partial(\Delta^1 \times \Delta^1 )) = V(\scrJ_{\alpha, \beta}) \cap V\la s_0s_1t_0t_1 \ra \ = \ \emptyset. 
\end{equation}
On the other hand, since
\[
\det \frac{\partial( \frp_\alpha , \frq_\beta )}{\partial( y_1, y_2)} =
\det \begin{pmatrix}
s_0s_1 & 0 \\
0 & t_0t_1
\end{pmatrix} \ = \ s_0s_1t_0t_1 
\]
it follows that the projection \( \spec{\scrA[y_1,y_2]/\scrJ_{\alpha, \beta}}\to \spec{\scrA} =X\times \Delta^1\times \Delta^1\) is smooth and surjects onto the complement of  \( X \times \partial(\Delta^1\times \Delta^1 )\). Furthermore, \(  \spec{\scrA[y_1,y_2]/\scrJ_{\alpha, \beta}}\) is irreducible and one can directly apply \cite[\S 7.2]{Ful-IT} to show that 
\(
 \text{div}(\frp_\alpha) \centerdot \text{div}(\frq_\beta) = [V(\scrJ_{\alpha, \beta}]
\)
is an irreducible cycle with multiplicity one, which is equidimensional and dominant over \( \Delta^1 \times \Delta^1 \). The second assertion now follows from \eqref{eq:inter1} and Remark \ref{rem:prism}.
 
Consider the projections
 \(\rho_1  \colon (X \times \bba^2) \times \Delta^1 \times  \Delta^1\times \bba^1 \to (X \times \bba^1) \times \Delta^1 \times \bba^1 \) 
 and
 \(\rho_2  \colon (X \times \bba^2) \times \Delta^1 \times  \Delta^1 \to (X \times \bba^1) \times \Delta^1 \) 
defined by 
\[
\begin{tikzcd}
   & \ar[dl, mapsto, "\rho_1"'] (x; y_1, y_2; s_0, s_1; t_0, t_1; T)  \ar[rd, mapsto, "\rho_2"] &   \\
(x; y_1; s_0, s_1; T)  &   &  (x; y_2; t_0, t_1).
\end{tikzcd}.
\]
The divisor
\(\rho_1^* H_{\alpha, \beta} \) intersects \( \rho_2^* D_\gamma\) properly and their cycle-theoretic intersection is equidimensional over \( \Delta^1 
\times \Delta^1\times \times \bba^1 \).  It follows directly from this observation and the proof of the first statement that
\[
\left(
\rho_1^* H_{\alpha, \beta} \centerdot\rho_2^* D_\gamma
\right)|_{T=0} = \pi_1^* D_{\alpha\beta} \centerdot\pi_2^* D_\gamma 
\]
and
\[
\left(
\rho_1^* H_{\alpha, \beta} \centerdot \rho_2^* D_\gamma
\right)|_{T=1} = \left( \pi_1^* D_{\alpha } + \pi_1^* D_{ \beta}\right)  \centerdot\pi_2^* D_\gamma  
=  \pi_1^* D_{\alpha }    \centerdot \pi_2^* D_\gamma  \ + \ 
 \pi_1^* D_{ \beta}  \centerdot \pi_2^* D_\gamma . 
\]
Applying the Eilberg-Zilber operator \( \psi_{1,1}^* \) to these identities, along with Remark~\ref{rem:prism}, concludes the proof of the third assertion. 

Recall that \( \la \frp_\alpha \ra + \la s_0s_1 \ra = \la 1 \ra = \la \frq_{\beta} \ra + \la t_0t_1 \ra \)  in \( \scrA[y_1,y_2]\), whenever \( \alpha, \beta \in R^\times \). In particular \( \la \frp_\alpha, \frq_{-\alpha} \ra = \la \frp_\alpha, \frg \ra \), where \( \frg = t_1^2(s_0^2 +s_0s_1y_1) + s_1^2 ( t_0^2 + t_0t_1y_2) \). Define
\( \hat \frg :=  t_1^2(s_0^2 +s_0s_1y_1) + s_1^2 ( t_0^2T^2 + t_0t_1 T y_2) \in \scrA[y_1, y_2, T] \) and
observe that \( Z:= \spec{\scrA[y_1,y_2,T]/\la \frp_\alpha, \hat \frg \ra} \subset (X \times \bba^2)\times \Delta^1 \times \Delta^1 \times \bba^1 \)
is smooth over \( \Delta^1 \times \Delta^1 \times \bba^1 \) and that \( Z \cap \left( X \times \bba^2 \times \Delta^1 \times \Delta^1 \times \{ 0 \}\right) = \emptyset \). It follows that \(G:=  \text{div}(\frp_\alpha) \centerdot \text{div}(\hat \frg) \) is an irreducible cycle of codimension \( 2 \) which is equidimensional and dominant over \( \Delta^1 \times \Delta^2 \times \bba^1 \) and satisfies \( G_{|T=0} =0 \) and \( G_{|T=1} = \text{div}(\frp_\alpha)\centerdot \text{div}(\frq_{-\alpha})\). Using previous arguments one concludes that \( \Gamma_{\alpha, -\alpha} = 0 \). 

When both \( \alpha\) and \( 1 - \alpha\) are in \(R^\times \) one can proceed in a similar fashion and define \( \hat \frg_1 := s_1^2\{ 1 + t_0 t_1T(y_2-2)\} - (t_1T)^2s^2_1 \alpha \), noting that  that \( \la \frp_\alpha, \frq_{1-\alpha} \ra = 
\la \frp_\alpha, \ s_1^2\{ 1 + t_0 t_1(y_2-2)\} - (t_1)^2s^2_1 \alpha \ra \). Define \( Z_1 :=  \spec{\scrA[y_1,y_2,T]/\la \frp_\alpha, \hat \frg_1 \ra} \subset (X \times \bba^2)\times \Delta^1 \times \Delta^1 \times \bba^1\) and \( G_1 :=  \text{div}(\frp_\alpha) \centerdot \text{div}(\hat \frg_1) \). Once again, 
 \(G \) is an irreducible cycle of codimension \( 2 \) which is equidimensional and dominant over \( \Delta^1 \times \Delta^1 \times \bba^1 \) and satisfies \( G_{1|T=0} =0 \) and \( G_{1|T=1} = \text{div}(\frp_\alpha)\centerdot \text{div}(\frq_{1-\alpha})\). Previous arguments imply that \( \Gamma_{\alpha, 1-\alpha} = 0 \). 

\end{proof}

\begin{corollary}
\label{cor:symbol}
The assignment 
\[
  \alpha\otimes \beta  \in R^\times\otimes R^\times\ \overset{\gamma}{\longmapsto}\  [ \Gamma_{\alpha,\beta}] =: \gamma_{\alpha,\beta} \in \shcg{2}{2}{X\times \bba^2} \cong \shcg{2}{2}{X} 
\] 
is a Steinberg symbol defining a natural transformation 
\[
\gamma 
\colon K_2^M(R)  \to \shcg{2}{2}{\spec{R}} 
\]
of functors in the category of regular \( \Bbbk\)-algebras. 
\end{corollary}

\appendix

\section{Nori's equidimensionality lemma (by Madhav Nori)}
\label{sec:LULU}

\noindent{\it The main theorem in this section is an unpublished result of  Madhav Nori, who kindly shared the proof and his insight with us. 
Its genesis goes back to a correspondence between Madhav and Ravi Rao that took place around $1988$.
}
\smallskip

All varieties and morphisms in this section are over \(  \spec{\Bbbk} \), for a base  field \( \Bbbk\).

\begin{theorem}[M. Nori]
\label{thm:equid}
Let  \(  G \) be an algebraic group of dimension \(  n,   \) 
and let \(  f_i \colon X_i \to G \), \(  i=1,\dots, k \), be dominant morphisms from irreducible varieties \(  X_i \) to \(  G. \) Define \(  X := X_1 \times \cdots \times X_k \) and \(  F\colon X \to G \) by \(  F(x_1, \ldots, x_k) = f_1(x_1)\cdots f_k(x_k). \) If \(  k\geq n \), then \(  F \) is an equidimensional morphism. 
\end{theorem}

\begin{remark}
\label{rem:sharp}
The bound is sharp. Let \(  f_i \colon X_i \to G \) be the blow-up of a point \(  p_i \in G, \) with exceptional fiber \(  E_i, i=1, \ldots k. \) Then the product of the exceptional fibers, which has dimension \(  k(n-1) \), is contained in \(  F^{-1}(p_1\cdots p_k). \) However, the general fiber has dimension \(  (k-1)n \), which is \(  < k(n-1) \) under the assumption \(  k<n. \)
\end{remark}

\begin{lemma}
\label{lem:flat}
If \(  f\colon X \to G \) is equidimensional and \(  c\colon W \to G \) is a morphism, then \(  (x,w)\mapsto f(x)c(w) \) defines and equidimensional  morphisms from \(  X\times W \) to \(  G. \)
\end{lemma}
\begin{proof}
The map \(  S \colon X \times W \to G\times W \) given by \(  S(x,w) = (f(x), w)  \) is equidimensional, the map \(  T \colon G\times W \to G\times W \) given by \(  T(g,w) = (g \cdot c(w), w)  \) is an isomorphism and the projection \( P \colon G \times W \to G\) is equidimensional. It follows that \(  P\circ S\circ T \) is equidimensional.
\end{proof}
\begin{notation}
\label{not:open}
By applying Grothendieck's theorem of generic flatness to \(  f_i \colon X_i \to G \), we obtain a nonempty open \(  U_i \subset G \) such that  \(  f_i \) restricts to a flat morphism \(  X_i':= f_i^{-1}(U_i) \to G.  \)  In particular, \(  f_i \colon X_i' \to G \) is equidimensional.
\end{notation}
The next result follows directly from the lemma.
\begin{corollary}
\label{cor:equids}
Under the hypothesis and notation from Theorem \ref{thm:equid}, the several restrictions of \(  F \) to
\(  A_i := X_1 \times X_2 \times \cdots \times X_i' \times \cdots \times X_k \ \to \ G, \, i=1, \dots, k,\) are all equidimensional morphisms.
\end{corollary}
\begin{notation}
\label{not:not2}
Denote \(  X' = \cup_{i=1}^k A_i \subset X  \) and  \(  W_i := X_i - X_i' , \, i=1, \ldots, k.\) In particular, \( W:= W_1\times \cdots \times W_k \) is the complement \(  W= X - X'. \)
\end{notation}
\begin{proof}(of Theorem \ref{thm:equid})
The irreducibility of \(  X_i \) shows that \(  \dim(X_i) \geq 1 + \dim(W_i). \) Summing over \(  i =1, \ldots, k \), we get
\(  \dim(X) \geq k + \dim(W) \). Subtract \(  n \) from both sides to get
\[
\operatorname{rel.dim.}(F) \geq (k-n) + \dim (W) \geq \dim (W), 
\]
since we are assuming \(  k\geq n. \)

Given \(  g\in G \), the fiber \(  F^{-1}(g) \) is the union of
\begin{enumerate}[(i)]
\item \(  (X' \cap F^{-1}(g) )  \), whose dimension is \(  \leq \operatorname{rel.dim.}(F) \) by the corollary, and
\item \( W \cap  F^{-1}(g)  \), whose dimension is \(  \leq  \dim(W) \leq \operatorname{rel.dim.}(F). \)
\end{enumerate}
It follows that \(  \dim(F^{-1}(g) ) \leq \operatorname{reldim.}(F). \) For non-empty fibers, the reverse inequality follows from the upper semi-continuity of the dimension of fibers. This concludes the proof.
\end{proof}

\section{Proof of technical lemmas}
\label{app:pf1}

\subsection{Proof of Lemma \ref{lem:L-smooth}}

%
We first show that a certain ``convex combination'' of matrices yields the various smooth maps that play a key role in our constructions.

\begin{lemma}
\label{lem:smooth}
Let \(   \mu_{n, r } \colon \mat{n}{n}^{\times r} \times \bba^r \to \mat{n}{n}   \) be the morphism that sends \(  R \)-valued points  \(  \mbA = (A_1| \cdots | A_r)  \in \mat{n}{n}^{\times r}(R)\) and \(  \mbt = (t_1, \ldots, t_r) \in \bba^r(R) \) to
\[
\mu_{n,  r}(\mbA, \mbt) :=  I + \sum_{s =1}^r t_s ( A_s - I), 
\]
where \(  I \) is the identity matrix.
\begin{enumerate}[a.]
\item \label{it:La} Define
\(
\hatmu_{n,r} \colon \mat{n}{n}^r \times  \bba^r  \to 
\mat{n}{n} \times   \bba^r  
\)
by\ \ \( \hatmu_{n,r}(\mbA, \mbt) = (\mu_{n,r}(\mbA, \mbt), \mbt) \).\ 
Then, the restriction of $\hatmu_{n,r}$ to $\mat{n}{n}^r \times \left( \bba^r - \{ 0 \} \right)$ is a smooth morphism.
In particular, the compositions
 \[
 \xymatrix{
&  \mat{n}{n}^r \times \left( \bba^r-\{ 0 \} \right) \ar[d]^-{\hatmu_{n,r}}
 \ar@/_1.6pc/@{-->}[dl]_-{\mu_{n,r}}  \ar@/^1.6pc/@{-->}[dr]  
 &    \\
 \mat{n}{n} & \mat{n}{n}\times \bba^r \ar[l]^-{pr_1} \ar[r]_-{pr_2}    
  &   \bba^r 
 }
 \]
are also smooth, where $pr_1$ and $pr_2$ are the indicated projections.
\item \label{it:Lb} Similarly, define 
\( 
\Psi_{n, r } \colon  \left( \bba^n -\{ 0 \} \right) \times \mat{n}{n}^{\times r} \times \bba^r \to \bba^n \times \bba^r
\)
on \(  R \)-valued points \(  (\mbu, \mbA, \mbt)  \) by \
\(
\Psi_{n, r }(\mbu, \mbA, \mbt)   := \left( \mbu \cdot \mu_{n, r}(\mbA, \mbt), \ \mbt    \right),
\)\ 
where \(  \mbu \cdot \mu_{n, r}(\mbA, \mbt)  \) denotes matrix multiplication when we see \(  \mbu  \) as an element in \(  \mat{1}{n}  \).  Then,  \(  \Psi_{n,r} \) is a smooth morphism. 
\end{enumerate}
\end{lemma}
\begin{proof}
Consider sets of variables \(  \mbx^1 = \{ x_{ij}^1 \}, \ldots, \mbx^r = \{ x_{ij}^r \}, \mby = \{ y_{ij} \}, \) with \(  1\leq i, j \leq n \),\ \(  \mbz = \{ z_1, \ldots, z_r \} \)
and denote \(  S:= \bbz[\mby, \mbz] \), so that \(\mat{n}{n} \times \bba^r = \mathsf{Spec}(S). \) 

It suffices to prove the first statement for the restrictions 
\(  \hatmu_{n,r{|\mat{n}{n}^r \times U_{\ell_0}}} \) where  \(  \bba^r-\{ 0 \} = \bigcup_{\ell=1}^r U_\ell\) is the open cover with
\(  U_\ell = \bba^r - V\la z_\ell\ra \cong \mathsf{Spec}(\bbz[\mbz  , \lambda ]/ \la z_\ell \lambda -1 \ra  \).
One can identify \(  \mat{n}{n}^{\times r} \times U_{\ell_0} \) with the graph of \(  \hatmu_{n,r{|\mat{n}{n}^r \times U_{\ell_0}}} \) and write
\[
\mat{n}{n}^{\times r} \times U_{\ell_0} \cong \mathsf{Spec}\left( S[\mbx^1, \ldots, \mbx^r,  \lambda ] / I \right),
\]
where
\(
I = \big\langle\ f_{ij}, \ \lambda\, z_{\ell_0} -1\ \mid \ 1\leq i, j \leq n\  \big\rangle,
\)
with
\(  f_{ij} := y_{ij} - \{\delta_{ij} +\sum_{\ell=1}^r z_\ell \, (x_{ij}^\ell - \delta_{ij}) \}  \).
In this way, the morphism \(  \hatmu_{n,r} \) is given by 
\(  \mathsf{Spec}\left( S[\mbx^1, \ldots, \mbx^r,   \lambda ] / I \right) \longrightarrow  \mathsf{Spec} (S).\)
Denote \(  g=\lambda z_{\ell_0} -1 \) and observe that
\[
\det \begin{pmatrix}
\frac{\partial f_{ij}}{\partial x_{rs}^{\ell_0}} & \frac{\partial g}{\partial x_{rs}^{\ell_0}} \\
 & \\
\frac{\partial f_{ij}}{\partial \lambda} & \frac{\partial g}{\partial \lambda} 
\end{pmatrix}
=
\det
\begin{pmatrix}
\delta_{ir}\delta_{js}\,  z_{\ell_0} & 0 \\
 & \\
 0 &  z_{\ell_0}
\end{pmatrix} \ = \ z_{\ell_0}^{n^2+1}. 
\]
Since \(  z_{\ell_0} \) is a unit in 
\(  S[\mbx^1, \ldots, \mbx^r,   \lambda ] / I \), it follows that \( \hatmu_{n,r} \) is smooth.

To prove the second assertion, it suffices to restrict  $\Psi_{n,r}$ to the open subsets
$U_s \times \mat{n}{n}^r \times \bba^r$, $s = 1, \ldots, n. $ First, consider  sets of variables $\mbu = \{ u_1, \ldots, u_n\}$, $\mbv = \{ v_1, \ldots, v_n\}$ so we can write $\bba^n \times \bba^r = \spec{R},$ with \( R = \bbz[\mbv, \mbz]\), and  let \( J \subset R[\mbu; \mbx^1, \ldots, \mbx^r; \lambda] \) be the ideal 
\[ 
J= \bigg\langle \left. v_k - \sum_{i=1}^n u_i \left\{ \delta_{ik} +\sum_{\ell=1}^r z_\ell ( x_{ik}^\ell -\delta_{ik} )\right\}, \ u_s \lambda -1 \ \right|  1\leq k \leq n \bigg\rangle.
\]
Denote 
\[
h_k :=v_k - \sum_{i=1}^n u_i \left\{ \delta_{ik} +\sum_{\ell=1}^r z_\ell ( x_{ik}^\ell -\delta_{ik})\right\}
\quad\quad  \text{ and } \quad \quad
g = u_s \lambda -1,
\]
and note that
\begin{equation}
    \label{eq:Jac}
    \frac{\partial h_k}{\partial x^a_{bc}} = - z_a u_b\, \delta_{kc} \quad \quad \text{and} \quad\quad  \frac{\partial h_k}{\partial u_a } = -\left( \delta_{ak} +  \sum_{\ell=1}^r z_\ell (x^\ell_{ik} - \delta_{ik}) \right).
\end{equation}

Let $\mathfrak{J}$ denote the ideal in $R[\mbu; \mbx^1, \ldots, \mbx^n, \lambda]/J$ generated by the maximal $(n+1)\times (n+1)$ minors of the matrix of partial derivatives of $h_1, \ldots, h_n, g$ with respect to the variables in  $\mbu\cup \mbx^1 \cup \cdots \cup \mbx^n \cup \{ \lambda\}$. Also, denote by $\scrI$ the ideal generated by $\{ z_1, \ldots , z_n \}\subset R$ and observe that, for  $1\leq k, j \leq n$, one has
\[
\det 
\begin{pmatrix}
\frac{\partial h_k}{\partial u_j} & \frac{\partial h_k}{\partial \lambda} \\
\frac{\partial g}{\partial u_j} & \frac{\partial g}{\partial \lambda} 
\end{pmatrix}
= \det
\begin{pmatrix}
-\delta_{k,j} + \scrI & 0  \\
\delta_{j,s}\lambda  & u_s
\end{pmatrix}
= (-1)^n u_s \ + \ \scrI,
\]
where we use the expression $ u= x+ \scrI$ to denote $u - x \in \scrI.$ 
Since $u_s$ is a unit, it suffices to show that $\scrI \subset \mathfrak{J}$. We first use (reverse) induction on the degree of monomials to show that $ z_1, \ldots, z_n \in \mathfrak{J} \cap R.$ Indeed, 
fix $s$ and choose $\mba = (a_1, \ldots, a_n)$ with $a_i \in \{ 1, \ldots, n \}$. Then, the determinant of the   matrix
\[
J(s;\mba) = \begin{pmatrix}
  \frac{\partial h_k}{\partial x^{a_j}_{s,j}} & \frac{\partial h_k}{\partial \lambda} \\
  \frac{\partial g}{\partial x^{a_j}_{s,j}} & \frac{\partial g}{\partial \lambda} 
\end{pmatrix}
=
\begin{pmatrix}
-z_{a_j}u_s \delta_{k,j} & 0 \\
 0 & u_s
\end{pmatrix}
\]
turns out to be $\det(J(s;\mba)) = (-1)^n u_s^{n+1} z_{a_1}\cdots z_{a_n}.$ Since $u_s$ is a unit, we conclude that 
$z_{a_1}\cdots z_{a_n} $ lies in $\mathfrak{J}$, and hence $\scrI^n \subset \mathfrak{J}.$ Now assume that $\scrI^k \subset \mathfrak{J}$ for all $k> d$ with $1 \leq d < n$. Given $\mbb = (b_1, \ldots, b_d)$
consider the matrix
\begin{align*}
M_d(s;\mbb) & = 
\left( 
\begin{array}{ccc|ccc|c}
\frac{\partial h_1}{\partial x^{b_1}_{s,1}} & \cdots & \frac{\partial h_1}{\partial x^{b_d}_{s,d}}  &
\frac{\partial h_1}{\partial u_{d+1}} & \cdots & \frac{\partial h_1}{\partial u_{n}} & \frac{\partial h_1}{\partial \lambda} \\
\vdots &  & \vdots & \vdots &  & \vdots & \vdots \\
\frac{\partial h_d}{\partial x^{b_1}_{s,1}} & \cdots & \frac{\partial h_d}{\partial x^{b_d}_{s,d}}  &
\frac{\partial h_d}{\partial u_{d+1}} & \cdots & \frac{\partial h_d}{\partial u_{n}} & \frac{\partial h_d}{\partial \lambda} \\
\frac{\partial h_{d+1}}{\partial x^{b_1}_{s,1}} & \cdots & \frac{\partial h_{d+1}}{\partial x^{b_d}_{s,d}}  &
\frac{\partial h_{d+1}}{\partial u_{d+1}} & \cdots & \frac{\partial h_{d+1}}{\partial u_{n}} & \frac{\partial h_{d+1}}{\partial \lambda} \\
\vdots &  & \vdots & \vdots &  & \vdots & \vdots \\
\frac{\partial h_n}{\partial x^{b_1}_{s,1}} & \cdots & \frac{\partial h_n}{\partial x^{b_d}_{s,d}}  &
\frac{\partial h_n}{\partial u_{d+1}} & \cdots & \frac{\partial h_n}{\partial u_{n}} & \frac{\partial h_n}{\partial \lambda} \\
\frac{\partial g}{\partial x^{b_1}_{s,1}} & \cdots & \frac{\partial g}{\partial x^{b_d}_{s,d}}  &
\frac{\partial g}{\partial u_{d+1}} & \cdots & \frac{\partial g}{\partial u_{n}} & \frac{\partial g}{\partial \lambda} \\
\end{array}
\right)  \\
& =
\left( 
\begin{array}{ccc|ccc|c}
-z_{b_1} u_s & \cdots & 0  &
\frac{\partial h_1}{\partial u_{d+1}} & \cdots & \frac{\partial h_1}{\partial u_{n}} & 0 \\
\vdots &  & \vdots & \vdots &  & \vdots & \vdots \\
0 & \cdots & -z_{b_d} u_s  &
\frac{\partial h_d}{\partial u_{d+1}} & \cdots & \frac{\partial h_d}{\partial u_{n}} & 0 \\ \hline
0 & \cdots & 0  &
\frac{\partial h_{d+1}}{\partial u_{d+1}} & \cdots & \frac{\partial h_{d+1}}{\partial u_{n}} & 0 \\ 
\vdots &  & \vdots & \vdots &  & \vdots & \vdots \\
0 & \cdots & 0  &
\frac{\partial h_n}{\partial u_{d+1}} & \cdots & \frac{\partial h_n}{\partial u_{n}} & 0 \\
0 & \cdots & 0  &
0 & \cdots & 0 & u_s 
\end{array}
\right) 
\end{align*}
Then 
\begin{align*}
\det \left( M_d(s;\mbb) \right) & = (-1)^d\, u_s^{d+1} z_{b_1}\cdots z_{b_d} \ 
\det \begin{pmatrix}
\frac{\partial h_{d+1}}{\partial u_{d+1}} & \cdots & \frac{\partial h_{d+1}}{\partial u_{n}} \\ 
\vdots &  & \vdots  \\
\frac{\partial h_n}{\partial u_{d+1}} & \cdots & \frac{\partial h_n}{\partial u_{n}} 
\end{pmatrix} \\
&
= (-1)^d\, u_s^{d+1}z_{b_1}\cdots z_{b_d}  \ 
\det \begin{pmatrix}
1+ \scrI & \cdots & \scrI \\ 
\vdots &  & \vdots  \\
\scrI & \cdots & 1+ \scrI 
\end{pmatrix} \\
& = 
(-1)^d\, u_s^{d+1}z_{b_1}\cdots z_{b_d} ( 1 + \scrI) = 
(-1)^d\, u_s^{d+1}z_{b_1}\cdots z_{b_d} + \scrI^{d+1}.   
\end{align*}
By definition, $\det \left( M_d(s;\mbb) \right) \in \mathfrak{J} $ and by induction $ \scrI^{d+1} \subset \mathfrak{J}$. Therefore, $z_{b_1}\cdots z_{b_d} \in \mathfrak{J}$ and hence $\scrI^d \subset \mathfrak{J}.$ This concludes the proof.
\end{proof}

\begin{proof}[Proof of Lemma \ref{lem:L-smooth}]

Let 
\(
 F_{n,r} \colon GL_n^{\times r} \xrightarrow{\, \cong\, } GL_n^{\times r}\) 
 \  and\  \(  p_0 \colon \Delta^n \xrightarrow{ \cong } \bba^n\) be the isomorphisms    defined 
  on \(  R \)-valued points $\mbA = (A_1|A_2|\cdots | A_r) \in GL_n^{\times}(R)$ and
$\mbt = (t_0, t_1, \ldots, t_n) \in \Delta^n(R)$ by
\begin{equation}
    \label{eq:isos}
    F_{n,r}( \uA ) := (A_1\mid A_1A_2\mid \cdots \mid A_1\cdots A_r) \quad\text{and}\quad p_0(\mbt) = (t_1, \ldots, t_n),
\end{equation}
 respectively.
It follows that the restriction under study can be written as the composition  in Figure \ref{fig:compo}
\begin{figure}[ht]
\(\xymatrix@C=3.6em{
GL_n^{\times r} \times \left( \Delta^r - \{ \mbe_0 \} \right) \ar[r]^{F_{n,r}\times p_0 \ \ }_{\cong}   \ar@/_1pc/[rd]_{\widehat L_{n,r}}
&  
GL_n^{\times r} \times \left( \bba^r - \{ 0 \} \right) \   \ar@{^{(}->}[r]^-{\jmath}   & \mat{n}{n}^r \times (\bba^r - \{ 0 \} ) \ar[d]^-{\hatmu_{n,r}} \ar[dl]_{\mu_{n,r} \times \mbp^{-1}_0\ }  \\
& \mat{n}{n}\times \Delta^r  & \ar[l]^{\bone\times \mbp^{-1}_0}  \mat{n}{n} \times \bba^r-\{ 0 \} ,
}
\)
\caption{}
\label{fig:compo}
\end{figure}
where \(  \jmath \) is the open inclusion. The first statement  now follows from Lemma \ref{lem:smooth}.

To prove the second statement,   note that the restriction \( F_{n,r} \colon SL_n^{\times r}  \xrightarrow{\cong} SL_n^{\times r}  \) is still an isomorphism. Hence,  it suffices to consider the composition below
\[
\begin{tikzcd}
SL_n^{\times r}  \times \Delta^r \arrow[r, "\mu"]  \arrow[rrr, "\mathsf{p}"',  bend right=12] & \mat{n}{n}\times \Delta^r \arrow[rr, "\pi_{n-1}\times 1"] & &   \mat{n}{(n-1)}\times \Delta^r,
\end{tikzcd}
\]
where \( \mu(\mbA,\mbt) = t_0I + t_1 A_1 + t_2 A_1A_2 + \cdots + t_r A_1 \cdots A_r\), and to show that the restriction of \( \mathsf{p}\) to 
\( SL_n^{\times r} \times \Delta^r - \{ \mbe_0\} \) is a smooth morphism. 

Consider sets of variables \(  \mbx^1 = \{ x_{ij}^1 \}, \ldots, \mbx^r = \{ x_{ij}^r \}\) with \(  1\leq i, j \leq n \),\ \( \mby = \{ y_{ij} \}, \) 
with \(  1\leq i \leq n,   2 \leq j \leq n\), 
 \(  \mbz = \{ z_1, \ldots, z_r \} \)
and denote \(  S:= \bbz[\mby, \mbz] \), so that \(\mat{n}{n-1} \times \bba^r = \mathsf{Spec}(S). \) 

For each \(k=1, \ldots, r\) set \( d_k := \det (\mbx^k) -1 \) and, given \( 1\leq i \leq n \),  define \( M_i(\mbx^k) = \det_{[n]-\{i \}, [n]-\{ 1 \} }(\mbx^k) \). 
Fix a multi-index \( I = (i_1,\ldots, i_r), i_k \in [n], k=1, \ldots, r \) and  \( \ell_\circ =1, \ldots, r \).
The scheme \( SL_n^{\times r } \times (\bba - \{ 0 \} ) \) can be covered by the basic open sets \( U_{I, \ell_\circ }\) where \(  z_{\ell_\circ}\) and \( M_{i_k}(\mbx^k), k=1, \ldots, r, \) are invertible. It follows that  the restriction \( \mathsf{p} \colon U_{I,\ell_\circ} \to \mat{n}{(n-1)}\times \Delta^r \cong 
\mat{n}{(n-1)}\times \bba^r \) is given by the homomorphism
\[
S \to S[\mbx^1, \ldots, \mbx^r,\lambda, \tau_1, \ldots, \tau_r]/\scrI_{I,\ell_\circ},
\]
where the ideal  
\(
\scrI_{I,\ell_\circ} := \la f_{i,j}, g,  h_k,  d_k \mid 1\leq i \leq n, \, 2\leq j \leq n,\, k=1, \ldots, r \ra,
\)
is generated by 
\begin{align*}
f_{ij} & := y_{ij} - \{\delta_{ij} +\sum_{\ell=1}^r z_\ell \, (x_{ij}^\ell - \delta_{ij}) \} \\
g & := \lambda z_{\ell_\circ} -1 \\
h_k & := \tau_k M_{i_k}(\mbx^k) - 1 \\
d_k & := \det(\mbx^k) -1 .
\end{align*}
It follows that (after appropriate ordering of indices)
\begin{align*}
& \det
\begin{pmatrix}
\frac{\partial f_{ij}}{\partial \mbx^{\ell_\circ}_{rs}}  &
\frac{\partial g }{\partial \mbx^{\ell_\circ}_{rs}} &
\frac{\partial h_k}{\partial \mbx^{\ell_\circ}_{rs}} &
\frac{\partial d_k }{\partial \mbx^{\ell_\circ}_{rs}} \\[0.25cm]
\frac{\partial f_{ij}}{\partial \lambda }  &
\frac{\partial g }{\partial  \lambda }&
\frac{\partial h_k}{\partial  \lambda } &
\frac{\partial d_k }{\partial  \lambda } \\[0.2cm]
\frac{\partial f_{ij}}{\partial \tau_a }  &
\frac{\partial g }{\partial \tau_a }&
\frac{\partial h_k}{\partial \tau_a } &
\frac{\partial d_k }{\partial \tau_a } \\[0.25cm]
\frac{\partial f_{ij}}{\partial \mbx^{a}_{i_a,1}}  &
\frac{\partial g }{\partial  \mbx^{a}_{i_a,1}}&
\frac{\partial h_k}{\partial  \mbx^{a}_{i_a,1} } &
\frac{\partial d_k }{\partial \mbx^{a}_{i_a,1}} \\[0.25cm]
\end{pmatrix} \\
&  = \ 
\det \begin{pmatrix}
\delta_{i,r}\delta_{s,r} \, z_{\ell_\circ} & 0 & * & * \\
0 & z_{\ell_\circ} & 0 & 0 \\
0 & 0 & \delta_{a,k} M_{i_a}(\mbx^a) & 0 \\
0 & 0 & 0 & (-1)^{i_a+1} \delta_{a,k} M_{i_a}(\mbx^a)
\end{pmatrix} \\
& =
(-1)^{r+|I|}z_{\ell_\circ}^{n^2-1}\prod_{a=1}^r M_{i_a}(\mbx^a).
\end{align*}
By definition of \( U_{I,\ell_\circ} \), this polynomial lies in \( \calo_{U_{I,\ell_\circ}}(U_{I,\ell_\circ})^\times \) and the lemma is proved. 
\end{proof}

\subsection{Proof of Proposition \ref{prop:prop}}
\label{subsec:pf_prop:prop}
Let \( \underline\mbA = (\mbA_0|\mbA_1|\cdots |\mbA_r) \) denote an element   in \( \scrF_{n,r} \) and let \( \frh_n \colon \scrF_n \times \bba^1 \to SL_n\) denote the contracting homotopy introduced in Definition \ref{def:contracting}. Then both \eqref{eq:restr} and assertion \ref{it:1}. in the proposition follow from the fact that    
\( \frh_n(\mbA, 0) = I \).

To prove assertion \ref{it:2}, let  \( \tilde \Gamma \colon X \times \scrF_{n,r}\times \Delta^{r+1} \to SL_n^{\times (r+1)} \) be  given by
\(
\tilde\Gamma =
\left(\tilde\gamma_1, \ldots, \tilde\gamma_r, \tilde\gamma_{r+1}
\right),
\)
where
\[
\tilde\gamma_i(x;\underline\mbA;\mbz) =
\begin{cases}
\mba_i(x;\, z_0,\ldots, z_{r-1}, z_{r}+z_{r+1})\, \frh_n(\mbA_i, \rho_{r+1}(\mbz)) &, \text{ if  }\  1\leq i \leq r \\
\frh_n(\mbA_0, \rho_{r+1}(\mbz)) &, \text{ if  }\   i = r+1.
\end{cases}
\]
A simple inspection shows that
\[
\partial_j \tilde\Gamma =
\begin{cases}
I &, \text{ if }\ 1\leq j \leq r-1, \\
\pi^*\mba &, \text{ if }\ j = r \\
\lambda_\mba&, \text{ if } \ j =r+1.
\end{cases}
\]

To prove assertion \ref{it:3}., let \( \gamma_i \colon X\times \Delta^{r+1} \to SL_n^{\times (r+1)}, 1\leq i \leq r+1, \) be the coordinate maps of a
 homotopy \( \Gamma \colon \mba\sim_h \mbb\). Define \( \tilde\Gamma \colon X \times \scrF_{n,r}\times \Delta^{r+1} \to SL_n^{\times (r+1)} \) with coordinate functions 
\[
\tilde\gamma_i(x;\underline\mbA; \mbz) := 
\begin{cases}
\gamma_i(x;\mbz) \, \frh_n(\mbA_i, \rho_r(\mbz) )  &, \text{ if } \ i=1, \ldots, r-1\\
\gamma_r(x;\mbz) \, \frh_n(\mbA_r, \rho_{r+1}(\mbz) ) &, \text{ if } i = r \\
\gamma_i(x;\mbz) \, \frh_n(\mbA_r, \rho_r(\mbz) ) \, \frh_n(\mbA_0;\rho(\mbz))&, \text{ if } i = r+1.
\end{cases}
\]

It follows from the definitions that \( \tilde\Gamma \colon \lambda_\mba \sim_h \lambda_\mbb\), i.e.
\[
\partial_j \tilde\Gamma =
\begin{cases}
I &, \text{ if }\ 1\leq j \leq r-1, \\
\lambda_\mba &, \text{ if }\ j = r \\
\lambda_\mbb&, \text{ if } \ j =r+1.
\end{cases}
\]
The assertion on flatness follows from Proposition \ref{prop:H}, Theorem \ref{thm:equid} and the definition of \( \tilde \Gamma\).

Finally, let \( \Gamma \colon X \times \Delta^{r+1}\to SL_n^{\times(r+1)}\) fill the horn \( (\underline I, \ldots, \underline I, \mba, - , \mbb ) \), and denote 
\( \mbc := \partial_r \Gamma \colon X\times \Delta^r \to SL_n^{\times r}\). Define \( \tilde \Gamma \colon X \times \scrF_{n,r} \times \Delta^{r+1} \to SL_n^{\times (r+1)} \) with coordinates defined by \( \tilde\gamma_i(x;\underline\mbA; \mbz) := \)
\[
\begin{cases}
\gamma_i(x;\mbz) \, \frh_n(\mbA_i, \rho_{r-1}(\mbz) ) \, \frh_n(\mbA_i, \rho_r(\mbz) ) \, \frh_n(\mbA_i, \rho_{r+1}(\mbz) )  &, \text{ if } \ 1\leq i \leq r-2\\
\gamma_{r-1}(x;\mbz) \, \frh_n(\mbA_{r-1}, \rho_{r}(\mbz) )\, \frh_n(\mbA_{r-1}, \rho_{r+1}(\mbz) ) &, \text{ if } i = r -1\\
\gamma_{r}(x;\mbz) \, \frh_n(\mbA_{r}, \rho_{r+1}(\mbz) )\, \frh_n(\mbA_{r-1}, \rho_{r-1}(\mbz) ) &, \text{ if } i = r \\
\gamma_{r+1}(x;\mbz) \, \frh_n(\mbA_r;\rho_{r-1}(\mbz))\, \frh_n(\mbA_r, \rho_r(\mbz) ) \, \frh_n(\mbA_0;\rho(\mbz))&, \text{ if } i = r+1.
\end{cases}
\]
It is easy to check that
\[
\partial_j \tilde\Gamma =
\begin{cases}
I &, \text{ if }\ 1\leq j \leq r-2, \\
\lambda_\mba &, \text{ if }\ j = r-1 \\
\lambda_\mbc &, \text{ if }\ j = r \\
\lambda_\mbb&, \text{ if } \ j =r+1.
\end{cases}
\]
Once again, the assertion on flatness follows from Proposition \ref{prop:H}, Theorem \ref{thm:equid} and the definition of \( \tilde \Gamma\).
\qed


\bibliographystyle{amsalpha}

\providecommand{\bysame}{\leavevmode\hbox to3em{\hrulefill}\thinspace}
\providecommand{\MR}{\relax\ifhmode\unskip\space\fi MR }
\providecommand{\MRhref}[2]{%
  \href{http://www.ams.org/mathscinet-getitem?mr=#1}{#2}
}
\providecommand{\href}[2]{#2}

\end{document}